\def\tr{\rm{tr}}
\def\Norm{\rm{Norm}}
\def\PSL{\rm{PSL}}
\def\SL{\rm{SL}}
\def\PGL{\rm{PGL}}

\def\SO{\rm{SO}}
\def\axis{\rm{axis}}

\documentstyle[12pt,psfig]{article}

\newtheorem{theorem}{Theorem}[section]    
                                           
\newtheorem{lemma}[theorem]{Lemma}         
\newtheorem{corollary}[theorem]{Corollary}

\makeatletter
\let \c@equation=\c@theorem
\makeatother

\title{Arithmeticity, Discreteness and Volume}
\author{F. W. Gehring \and  C. Maclachlan \and G. J. Martin \and A. W. Reid 
\thanks{Research supported in part by grants from the U. S. National Science
Foundation,  the N.Z. Foundation of Research,  Science and Technology, the
Australian Research Council, the U. K. Royal Society and the U.K. Scientific and
Engineering Research Council. We wish also to thank the University of Texas at
Austin, the University of Auckland and the MSRI at Berkeley (NSF Grant DMS-90222140) for
their support during part of this work.  We are grateful to  D. J. Lewis who gave us a
proof for Theorem 5.12 and to K. N. Jones who computed the co-volumes that appear in
Table 11.}}

\begin{document}           

\maketitle

\begin{abstract}
We give an arithmetic criterion which is sufficient to imply the discreteness of various
two-generator subgroups of $PSL(2,{\bf C})$.  We then examine certain two-generator
groups which arise as extremals in various geometric problems in the theory of Kleinian
groups, in particular those encountered in efforts to determine the smallest co-volume,
the Margulis constant and the minimal distance between elliptic axes.  We establish the
discreteness and arithmeticity of a number of these extremal groups,  the associated
minimal volume arithmetic group in the commensurability class and we study whether or
not the axis of a generator is simple. 
\end{abstract}

\section{Introduction}

In the papers \cite{GM1},  \cite{GM3}, \cite{GM4}, \cite{GM5}, \cite{GM7}
concerning geometric constraints on Kleinian groups, we identify many two-generator
subgroups of $PSL(2,{\bf C})$ which are candidates for various extremal problems. 
These include a family of extreme groups in which the minimum distance between
non--intersecting elliptic axes is realized, the minimum co-volume Kleinian group, the
Kleinian group in which the Margulis constant is realized, and a family of discrete
groups in which the distance between the fixed points of tetrahedral, octahedral and
icosahedral subgroups is also minimized; see also \cite{DM}. Not surprisingly all of
these problems are closely related.

A principal obstruction to solving some of these questions was in establishing the
discreteness of proposed candidates for the extremals.  Often these groups were either
two-generator groups or contained two-generator subgroups of finite index. In several
cases we established the discreteness of these two-generator groups by an ad hoc method;
see \cite{GM3}.  A main result of this paper gives a unified method, derived from
elementary algebraic number theory, which deals with most of these extremal situations. 
As a consequence we find that many of the extremals are subgroups of arithmetic Kleinian
or Fuchsian groups. These results are derived from a classical result in number theory
whose application to proving groups discrete in $PSL(2,{\bf C})$ appears to be new. When
one considers how difficult it is in general to prove that a two-generator subgroup of
$PSL(2,{\bf C})$ is discrete, the criterion considered here seems surprisingly
straightforward and elegant.

Most of the extremal examples mentioned above arise from a disk covering
procedure based on complex iteration and a special semigroup of polynomials.  
This method is ideally suited for the discreteness criteria which we shall
establish here.  Hence we first outline this procedure in \S 2 in order to motivate
the criteria which we then give and discuss at the end of \S 2; \S 3, \S 4 and 
\S 5 contain background material on arithmetic groups and the proofs for the
results given in \S 2.

In \S 6 we apply these criteria to a family of two-generator groups $\langle f,g
\rangle$ where $f$ and $g$ are elliptics of orders $n$ and $2$, $n=3,4,5,6,7$. 
These groups resulted from the disk covering method discussed above while studying
the possible distances between non-intersecting elliptic axes in a discrete group,
possible minimum co-volume Kleinian groups, the Margulis constant and the length
structure of closed geodesics.

The groups $\langle f,g \rangle$ in \S 6 are subgroups of arithmetic groups.  In \S 7
and \S 8 we first compare the minimal co-volume of the arithmetic groups in which the
groups $\langle f,g \rangle$ embed with the distance between the axes of the generators
$f$ and $g$.  Using a  computer program developed by K. Jones and the last author it is
possible to determine which of the groups $\langle f,g \rangle$ have finite co-volume
and hence are themselves arithmetic.

Finally motivated by the small volume problem for orbifolds, we derive in \S 9 and then
apply in \S 10 arithmetic and geometric criteria to determine which of the groups
$\langle f,g \rangle$ of \S 6 have $f$ as a simple elliptic generator.

\section{Disk covering method and discreteness}

Let $G=\langle f,g \rangle$ be a two-generator subgroup of $PSL(2,{\bf C})$.  We
associate to each such group three complex numbers called the {\em parameters} of
$G$,  \begin{equation}
{\rm par}(G) = (\gamma(f,g),\beta(f),\beta(g))\in {\bf C}^3,
\end{equation}
where
\[\beta(f)={\tr}^2(f)-4,\\\ \beta(g)={\tr}^2(g)-4,\\\
\gamma(f,g)={\tr}([f,g])-2\]
and $[f,g]=fgf^{-1}g^{-1}$ is the multiplicative commutator.  These three
complex numbers determine $G$ uniquely up to conjugacy whenever $\gamma(f,g)\neq
0$; recall $\gamma(f,g)=0$ implies the existence of a common fixed point. 
Conversely, every such triple determines a two-generator subgroup of $PSL(2,{\bf C})$.

When $\gamma\neq 0,\beta$, there is a natural projection from ${\bf C}^3$
to ${\bf C}^2$, given by
\[  (\gamma,\beta,\beta^\prime) \rightarrow (\gamma,\beta,-4), \]
which preserves discreteness; see \cite{GM3}.  That is, if
$(\gamma,\beta,\beta^\prime)$ are the parameters for a discrete two-generator group,
then so are $(\gamma,\beta,-4)$.  Hence in these circumstances one can always replace
one of the generators in a discrete two-generator group by an element of order two
without altering the commutator or discreteness.   No such result holds for groups with
more than two-generators \cite{Ca}.

We shall study the pairs of parameters $(\gamma(f,g),\beta(f))$ in ${\bf C}^2$
which can correspond to discrete groups $\langle f,g \rangle$.  For the purposes of
this discussion we specialize to the case where $f$ is a primitive elliptic
element of order $n\geq 3$. Then \[\beta =\beta(f)=-4\sin^2(\pi/n)\]
and for each $n\geq 3$, we want to describe the set $E_n$ of values in the
complex plane which may be assumed by the commutator parameter $\gamma(f,g)$ in
the case when $<f,g>$ is discrete. This set cannot have 0 as a limit point \cite{GM3}.

Next if
\begin{equation}
h = g \circ f^{m_1} \circ g^{-1} \circ f^{m_2} \circ g \ldots f^{m_n} \circ
g^{\epsilon}
\end{equation}
where $\epsilon=(-1)^{n}$, then it follows from \cite{GM3} that
\begin{equation}
\gamma(f,h) =p(\gamma(f,g),\beta(f))
\end{equation}
where $p$ is a polynomial in the two variables with integer coefficients.
Equation (2.3) shows that for such words $h$, the value of $\gamma(f,h)$ depends
only on $\gamma(f,g)$ and $\beta(f)$ and {\em not} on $\beta(g)$. In particular,
the family of such polynomials $p$ is closed under composition in the first
variable and so forms a polynomial semigroup {\bf P}. Thus we may study the
iterates of a specific polynomial and, in this way, generate a sequence of
commutator traces which cannot accumulate at $0$.

For example if $h=gfg^{-1}fg$ and $\gamma = \gamma(f,g)$, then
\[p(\gamma,\beta)=\gamma(1+\beta - \gamma)^2.\]
Hence if $f$ is elliptic of order $6$, then $\beta=-1$ and
$p(\gamma,\beta)=\gamma^3$.
In this case if we set $\gamma_j=\gamma(f,g_j)$ where
\[ g_{j+1}=g_jfg_j^{-1}fg_j \;\;\; {\rm and} \;\;\; g_1=g,\] 
then $\gamma_{j+1} = \gamma^{3^j} \rightarrow 0$ as $j \rightarrow \infty$
if $0<|\gamma|<1$, a contradiction.  Thus $\gamma=0$ or $|\gamma|\geq1$ and
we obtain an analogue for elliptics of order 6 of the classical inequality of
Shimizu--Leutbecher which yields the same conclusion whenever $f$ is parabolic
\cite{Mas}.  In particular, 
\[E_6 \subset {\bf C} \setminus D_6\]
where $D_6 = \{z: 0<|z|<|\}$.

For each value of $n \geq 3$ we can argue as above using iteration with an
appropriate polynomial in ${\bf P}$ to obtain an open punctured disk $D_n$
about $0$ which does not contain the parameter $\gamma(f,g)$.  When $n\geq 7$, the
existence of such a neighborhood follows from J\o rgensen's inequality. Next if
$p(\gamma)=p(\gamma,\beta)$ is the polynomial in ${\bf P}$ which corresponds to
the element $h$ in (2.2), then the subgroup $\langle f,h \rangle$ cannot be discrete if
\[p(\gamma(f,g))=p(\gamma(f,g),\beta(f)) = \gamma(f,h) \in D_n.\]
Hence $p^{-1}(D_n)$ is an excluded region for $\gamma(f,g)$ and we conclude that
\[E_n \subset {\bf C} \setminus \{ \cup p^{-1}(D_n):p \in {\bf P}\}.\]
In this way we can cover regions of the complement of the set $E_n$ of admissible
values for the commutator parameter of groups with a generator of order $n$.

The program for describing these regions was begun in \cite{GM3} and was
significantly extended in \cite{GM6}.  Obtaining various geometric constraints for
such groups amounts to showing that a certain region $\Omega$ lies in the
complement of $E_n$ except possibly for some {\em exceptional} values. 

The result of this disk covering argument is illustrated below for the case where
$n=3$. All exceptional values within the region covered by the disks are illustrated and
{\em every} value turns out to correspond to a subgroup of an arithmetic group.  Many
are in fact two-generator arithmetic groups themselves.  That is, they are additionally
of finite co-volume.  
\bigskip

\centerline{\bf Possible values for commutator parameter when n=3}
\nopagebreak
\psfig{file=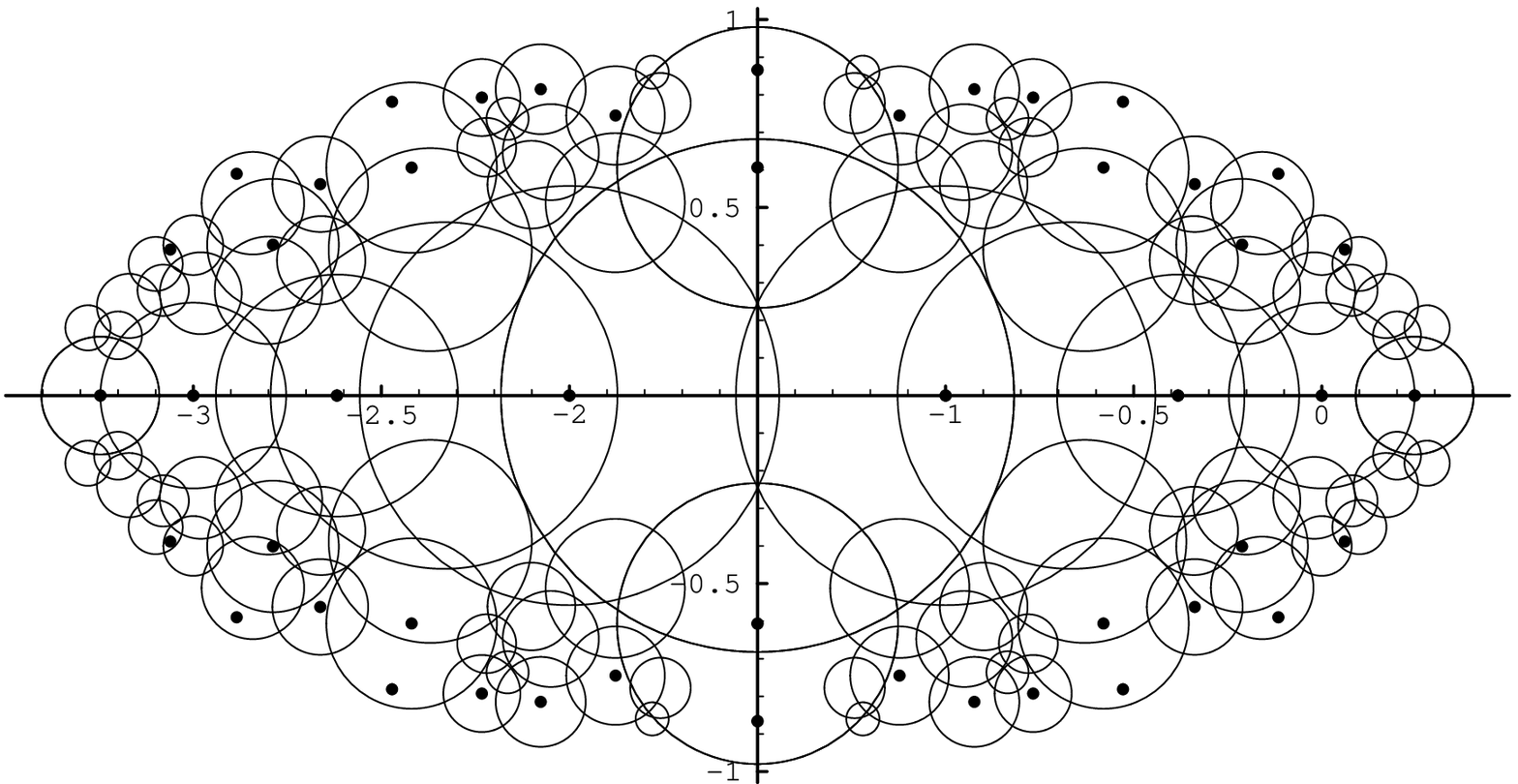,width=5.42in}

\bigskip

The closure of the hyperbolic line joining the two fixed points of a non-parabolic
element $f$ is called the {\it axis} of $f$, denoted by ${\axis(f)}$.  The following
formula yields collaring theorems from which volume estimates may be obtained; see
\cite{GM4} and \cite{GM8}. 

\begin{lemma} 
If $f$ and $g$ are non-parabolic, then the
hyperbolic distance $\delta(f,g)$ between the axes of $f$ and $g$ is given by   
\begin{equation}
\cosh(2\delta(f,g)) =
\left|\frac{4\gamma(f,g)}{\beta(f)\beta(g)}+1\right|+
\left|\frac{4\gamma(f,g)}{\beta(f)\beta(g)}\right|.
\end{equation}
\end{lemma}

The part of $E_n$ in $\Omega$ is contained in the union $X_n$ of the preimages of the
punctures of $D_n$ taken over all the polynomials $p$ in ${\bf P}$.  Moreover
$E_n$ is a proper subset of $X_n$ since the preimage of a puncture under one
polynomial $p$ in ${\bf P}$ can often be mapped into $D_n$ by some other
polynomial $q$ in ${\bf P}$, and hence will not be a point in $E_n$.  At the 
end of such a covering argument one is left with a few points of $X_n$ inside the
given region $\Omega$ which {\em may} correspond to discrete groups.  One must
then analyze these exceptions and prove that either they are not discrete or, if
they are discrete, use knowledge of these exceptional groups to compute the
invariants being studied.

\bigskip

It is often not easy to decide whether a given point $\gamma$ in $X_n$ belongs to
$E_n$, i.e. whether $\gamma$ corresponds to discrete group with a generator of
order $n$. However in our situation there exists a certain  polynomial with
coefficients in ${\bf Z}[\beta]$ for which $\gamma$ is a root. Thus it is  natural
to look for arithmetic conditions which guarantee discreteness.  The fortunate
fact is that this works well for the small extremals, the general philosophy being
that ``small'' implies ``highly symmetric'' which in turn implies ``arithmetic''.

Under the assumptions of discreteness and finite co-volume on subgroups of 
$PSL(2,{\bf R})$ and $PSL(2,{\bf C})$, the results of \cite{Ta} and \cite{MR}
give conditions on the traces of elements in the subgroup which 
characterize, among Fuchsian and Kleinian groups, those which are 
arithmetic; see \S 4. Dropping the assumption of finite co-volume 
characterizes those which are subgroups of arithmetic Fuchsian or 
Kleinian groups. More importantly, the proofs in \cite{Ta} and \cite{MR} 
make no use of the discreteness assumption; cf. \cite{B}. Thus the 
arithmetic conditions on the traces can be used to deduce discreteness.

In this paper refinements of these arithmetic conditions which will 
guarantee discreteness are obtained. These are particularly applicable to 
two-generator groups such as those corresponding to the points 
$\gamma$ in $X_n$ described above.

For example, suppose that  $G = \langle f , g \rangle$ where $f$ is a primitive elliptic
element of order $n \geq 3$ and $g$ is elliptic of order 2.  Then
\[\beta = \beta(f) = -4 \sin^2(\pi/n) \] 
is an algebraic integer in the field ${\bf Q}(\cos(2\pi/n))$, which is totally
real with $\phi(n) / 2$ places. Let $R_n$ denote the ring of integers in 
${\bf Q}( \cos(2\pi/n))$ and note that $R_n = {\bf Z}[\beta]$. The Galois
conjugates of $\beta$, $\sigma(\beta) = -4 \sin^2(m\pi/n)$ where $(m,n) = 1$,
lie in the interval $(-4 , 0)$ with $\sigma(\beta) \leq \beta$. Recall
that $G$ is determined up to conjugacy by $\beta$ and
$\gamma(f,g)$. With this we prove the following criteria for
discreteness. Recall that, for a number field $K$ of degree $n$ over
$\bf Q$, the $n$ Galois embeddings $\sigma : K \rightarrow {\bf C}$
give rise to valuations on $K$ which fall into equivalence classes -
the {\em places} - modulo the action of complex conjugation.

\begin{theorem}
Let $G = \langle f, g \rangle$ be a subgroup of $PSL(2,{\bf C})$ 
with $f$ a primitive elliptic element of order $n \geq 3$, $g$ 
an elliptic of order $2$  and $\gamma(f,g) \neq 0, \beta(f)$. Then 
$G$ is a subgroup of a discrete arithmetic group if 
\begin{enumerate}
\item ${\bf Q}(\gamma(f,g),\beta(f))$ has at most one complex place;
\item $\gamma(f,g)$ is a root of a monic polynomial $p(z) \in {\bf Z}[\beta][z]$;
\item if $\gamma(f,g)$ and $\bar{\gamma}(f,g)$ are not real, then all other roots of 
$p(z)$ are real and lie in the interval $(\beta,0)$;
\item if $\gamma(f,g)$ is real, then all other roots of $p(z)$ are real and 
lie in the interval $(\beta,0)$.
\end{enumerate}
\end{theorem}

The requirement that the field ${\bf Q}(\gamma(f,g),\beta(f))$ has at most one 
complex place can, for $\beta(f)$ as in Theorem 2.6, also be described  in terms of
the polynomial $p(z)$. If in addition $n = 3,4$ or $6$, then $\beta \in {\bf Z}$ and
the criteria admit a greatly simplified  description. See Theorem 5.14 below.

\section{A Criterion for Discreteness}

In this section we prove a result that guarantees discreteness under certain
conditions on the invariant trace field.  Our result is related to a classical
result in number theory and does not seem to be well known in the area of
Kleinian groups. Indeed it is a generalization of the prototype result that yields
discreteness for arithmetic subgroups of $PSL(2,{\bf C})$.  As we point out in
\cite{JR} this theorem has a proof in the language of arithmetic groups which
makes application more routine for our purposes; see also the discussion in \S 4.
However we give  here an elementary proof which illuminates the connection with
the following well known facts from number theory.

Recall that if $p$ is a monic irreducible polynomial over $\bf Z$ of degree
$n$ with roots $\alpha_1, \ldots ,\alpha_n$, then
\begin{equation}
p(z) = \prod_{i=1}^n (z- \alpha_i) = z^n - s_1z^{n-1} + \ldots
(-1)^ks_kz^{n-k} + \ldots
(-1)^ns_n
\end{equation}
where $s_i$ is the $i^{\rm th}$ symmetric polynomial in  $\alpha_1, \ldots
,\alpha_n$.  As a consequence we deduce the following easy lemma.

\begin{lemma}
There are only finitely many algebraic integers $z$ of bounded degree such
that $z$ and all Galois conjugates of $z$ are bounded.
\end{lemma}
{\bf Proof.}
Let $z$ have degree $n$ over ${\bf Q}$. If $z$ and its Galois conjugates
are bounded, then the coefficients of the irreducible polynomial $p$ of $z$ are
symmetric polynomials in the roots of $p$ and hence are bounded. Therefore only
finitely many integers can arise as coefficients of such a polynomial. $\Box$

\bigskip
Next we recall some notation from \cite{NR}.  Let $G$ be a finitely
generated subgroup of $PSL(2,{\bf C})$.  The {\em trace field} of $G$ is the field
generated over ${\bf Q}$ by the set 
\[{\tr}(G) =\{\pm {\tr}(g) : g \in G\}. \]
Since  $G$ is finitely generated, the subgroup $G^{(2)} = 
\langle g^2:g\in G \rangle$ is a normal subgroup of finite index with quotient
group a finite abelian $2$-group. Following \cite{NR} we call
\[  kG = {\bf Q}({\tr}(G^{(2)})) \]
the {\it invariant trace-field} of $G$.  For any finite index subgroup
$G_1$ of a nonelementary group $G$ one can show that  
${\bf Q}({\tr}(G^{(2)})) \subset {\bf Q}({\tr}(G_1))$; in \cite{R2} it is shown
that $kG$ is an invariant of the commensurability class. 

Throughout the paper, we will use 
\begin{equation}
{\bf c} : {\bf C} \rightarrow {\bf C}
\end{equation}
to denote the 
complex conjugation map. With this notation we can establish the following criterion for
discreteness.

\begin{theorem}
Let $G$ be a finitely generated subgroup of $PSL(2,{\bf C})$ such that
\begin{enumerate}
\item $G^{(2)}$ contains elements $g_1$ and $g_2$ which have no common
fixed point;
\item ${\tr}(G)$ consists of algebraic integers;
\item for each embedding $\sigma : kG \rightarrow {\bf C}$ such that 
$\sigma \neq id ~ or ~ {\bf c}$, the set $\{ \sigma(\tr(f)): f \in
G^{(2)}\}$ is bounded.
\end{enumerate}
Then $G$ is discrete.
\end{theorem}
{\bf Proof.}
Note that since $G$ is finitely generated, so is $G^{(2)}$ and so all traces 
in $G^{(2)}$ are obtained from integral polynomials in a finite number 
of traces. Thus $kG$ is a finite extension of ${\bf Q}$.

It suffices to prove that the finite index subgroup $G^{(2)}$ is discrete.
Suppose that this is not the case and let $f_n$ be a sequence of distinct elements
converging to the identity in $G^{(2)}$.  If  $z_n = {\tr}(f_n)$ and $z_{n,i} =
{\tr}([f_n,g_i])$, then 
\[ \beta(f_n)=z_{n}^{2}-4 \rightarrow 0 \\\ {\rm and} \\\
\gamma(f_n,g_i)=z_{n,i}-2 \rightarrow 0 \]
for $i=1,2$ as $n \rightarrow \infty$.  Hence we may assume that $|z_n| <
K$ for some fixed constant $K$. Next by condition 3., $|\sigma(z_n)| <
K_{\sigma}$ for each embedding $\sigma \neq id~ or ~{\bf c}$ of $kG$, where  $K_{\sigma}$
is a constant which depends only on $\sigma$.

Let $R = \max \{K, K_{\sigma}\}$ where $\sigma$ ranges over all
embeddings  $\sigma \neq id~or~{\bf c}$ of $kG$. Then the algebraic integers $z_n$ are of
bounded degree and they and all their Galois conjugates are bounded in absolute
value by $R$. By Lemma 3.2, the $z_n$ assume only finitely many values.  Thus for
large $n$, $\beta(f_n) = 0$ and $f_n$ is parabolic with a single fixed point
$w_n$.

Next we can apply the above argument to the algebraic integers $z_{n,i}$ to
conclude that $\gamma(f_n,g_i) = 0$  for $i=1,2$ and large $n$.  This then implies
that $g_1$ and $g_2$ each have $w_n$ as a common fixed point for large $n$
contradicting condition 1. $\Box$

\section{Arithmetic Kleinian Groups}

In this section we first recall some terminology about quaternion algebras
and then discuss Theorem 3.4 in the language of these algebras and arithmetic
groups. In particular we give a version of Theorem 3.4 which is readily
applicable, especially in the context of two-generator groups; see \S 5.

We begin with some facts about quaternion algebras; see \cite{V} for
details.  Let $k$ be a number field, let $\nu$ be a place of $k$, i.e. an
equivalence class of valuations on $k$ and denote by $k_{\nu}$ the completion of
$k$ at $\nu$. If $B$ is a quaternion algebra over $k$, we say that $B$ is {\it
ramified} at  $\nu$ if $B \otimes_k k_{\nu}$ is a division algebra of
quaternions. Otherwise $\nu$ is unramified.

In case $\nu$ is a place associated to a real embedding of $k$, $B$ is
ramified if and only if $B \otimes_k k_{\nu} \cong {\cal H}$, where $\cal H$ is the
Hamiltonian division algebra of quaternions.

It is straightforward to check whether a quaternion algebra is ramified at a
real place.  Recall, following \cite{V}, the {\it Hilbert symbol} $({{a,b}\over
k})$ corresponding to the quaternion algebra $\{1,i,j,ij\}$ defined over the field
$k$ with $i^2=a, j^2=b$ and $ij=-ji$.

\begin{lemma}
Let $B = ({{a,b}\over k})$. Then $B$ is ramified at a real place $\nu$
corresponding to the real embedding $\sigma$ if and only if $\sigma (a)$ and
$\sigma (b)$ are negative. 
\end{lemma}
{\bf Proof.}
By definition, $B$ is ramified at $\nu$ whenever $B \otimes_k k_{\nu} \cong 
{\cal H}$. This tensor product is isomorphic to \[ \biggl({{\sigma (a),\sigma (b)}\over
\sigma(k)} \biggr) \otimes_{\sigma (k)}{\bf R}\]
or simply $({{\sigma (a),\sigma (b)}\over {\bf R}})$.  This, in turn, is isomorphic
to $\cal H$ exactly when $\sigma(a)$ and $\sigma(b)$ are negative since one can
remove squares without affecting the isomorphism class of the quaternion algebra;
cf. \cite{V}. $\Box$

\bigskip
We now give the definition of an {\it arithmetic} Kleinian group. Let
$k$ be a number field with one complex place and $A$ a quaternion
algebra over $k$ ramified at all real places. Next let $\rho$ be an
embedding of $A$ into $M(2,{\bf C})$, $\cal O$ an order of $A$ and
${\cal O}^1$ the elements of norm 1 in $\cal O$.  Then $\rho({\cal
O}^1)$ is a discrete subgroup of $SL(2,{\bf C})$ and its projection to
$PSL(2,{\bf C})$, $P\rho({\cal O}^1)$, is an arithmetic Kleinian group.
The commensurability classes of arithmetic Kleinian groups
are obtained by considering all such $P\rho({\cal O}^1)$.

Arithmetic Fuchsian groups arise in a similar manner. In this case
the number field is totally real and the algebra ramified at all real places
except the identity.

In \cite{MR} (resp. \cite{Ta}), it is shown that two arithmetic Kleinian groups (resp.
arithmetic Fuchsian groups) are commensurable up to conjugacy if and only if their
invariant quaternion algebras are isomorphic; see also \cite{Bo}.

To state the characterization theorems in \cite{MR} and \cite{Ta} referred to
earlier, we first define, following Bass \cite{B},  
$$AG = \{ \sum a_ig_i: a_i \in {\bf Q}({\tr}(G)), \quad g_i \in G \}$$ 
for any finitely generated non-elementary subgroup of $SL(2,{\bf C})$.  Then $AG$ is a
quaternion algebra over ${\bf Q}({\tr}(G))$. Recall that $kG = {\bf Q}({\tr}(G^{(2)}))$
and $AG^{(2)}$ is then the invariant quaternion algebra since, for a finitely generated
non-elementary subgroup $G$ of $PSL(2,{\bf C})$, the pair $(AG^{(2)},kG)$ is an
invariant of the commensurability class of $G$. See \cite{NR}. Additionally, if
${\tr}(G)$ consists of algebraic integers, then as in \cite{B} 
$${\cal O}G = \{ \sum a_ig_i \mid a_i\in R_{{\bf Q}({\tr}(G))}, \quad g_i \in G
\}$$  
is an order in $AG$, where $R_{{\bf Q}({\tr}(G))}$ denotes the ring of integers
in ${\bf Q}({\tr}(G))$.

\begin{theorem}
Let $G$ be a Kleinian (resp. Fuchsian) group of finite co-volume. Then 
$G$ is arithmetic if and only if the following conditions are satisfied:
\begin{enumerate}
\item $kG$ is an algebraic number field;
\item ${\em tr}(G)$ consists of algebraic integers;
\item for every ${\bf Q}$-isomorphism $\sigma : kG \rightarrow {\bf C}$ 
such that $\sigma \neq id,{\bf c}$, $\sigma({\em tr}(G^{(2)}))$ is 
bounded in ${\bf C}$.
\end{enumerate}
\end{theorem}

\noindent{\bf Sketch of Proof}.  It is a simple matter to show that arithmetic
Kleinian and Fuchsian groups have these properties. Note that if
$\gamma \in {\cal H}^1 \cong S^3$ then $\tr(\gamma) \in [-2, 2]$.

Since $G$ has finite co-volume, it is finitely generated and non-elementary.
Only these two properties are used in the proof, except at the final stage.
Condition 1. actually follows from 2. and the fact that $G$ is finitely 
generated. We thus obtain a quaternion algebra $AG^{(2)}$ over $kG$ with $AG^{(2)}
\subset M(2,{\bf C})$. Furthermore,  each embedding $\sigma : kG \rightarrow {\bf C}$
extends to an embedding of $AG^{(2)}$ in $M(2,{\bf C})$. 

If $\gamma \in SL(2,{\bf C})$ has eigenvalues $u, u^{-1}$ whose
absolute values are not one, then it follows that $\{{\tr}(\gamma^m): m \in {\bf Z}\}$
is unbounded since  
\[|{\tr}(\gamma^m)| \geq ||u|^m -|u|^{-m}|.\]  
Thus in the case at hand condition 3. implies that $kG$ has exactly one complex place
if $kG$ is not real and is totally real if $kG$ is real. Also with $\gamma$
as above, $\tr(\gamma) \in [-2,2]$ and condition 3. also implies that
$AG^{(2)}$ is ramified at all real places of $kG$ $\neq id$. Finally
condition 2. shows that ${\cal O}G^{(2)}$ is an order of $AG^{(2)}$
and, of course, $G^{(2)} \subset ({\cal O}G^{(2)})^1$, which is a
discrete arithmetic Kleinian or Fuchsian group. Since $G$ has finite
co-volume, it will be commensurable with $({\cal O}G^{(2)})^1$ and so
be arithmetic. $\Box$

\bigskip

In practice, direct application of Theorem 4.2 (or Theorem 3.4) is hard, 
the problem being to establish the boundedness of the traces at real
embeddings. However, as a corollary of the proof of Theorem 4.2, we
obtain the following more useful method for proving groups discrete.

\begin{theorem}
Let $G$ be a finitely generated non-elementary subgroup of $PSL(2,{\bf C})$
such that 
\begin{enumerate}
\item $kG$ has exactly one complex place or is totally real;
\item $\tr(G)$ consists of algebraic integers;
\item $AG^{(2)}$ is ramified at all places of $kG$, $\neq \{id,{\bf c}\}$.
\end{enumerate}
Then $G$ is a subgroup of an arithmetic Kleinian or Fuchsian group.
\end{theorem}

If $g_1,g_2 \in G^{(2)}$ have no common fixed point and $g_1$ is not
parabolic, the proof of Theorem 4.2  shows that $AG^{(2)}$ is spanned by
$1,g_1,g_2,g_1g_2$ over $kG$; see \cite{HLM}, \cite{MR} and \cite{Ta}. A basis in
standard form can then be obtained yielding the Hilbert symbol $$\left( {\beta(g_1) ,
\gamma(g_1,g_2) \over kG} \right).$$ Now if $f,g \in G$ are a pair of non-commuting
elements with $f$ not parabolic and $f,g$ not of order 2, then $$\beta(f^2) =
(\beta(f) + 4)\beta(f), \quad \gamma(f^2,g^2) = (\beta(f) + 4)
(\beta(g) + 4)\gamma(f,g)$$ so that
\begin{equation}
AG^{(2)} \cong \left( {(\beta(f) + 4)\beta(f), (\beta(f) + 4)(\beta(g) + 4)
\gamma(f,g) \over kG} \right).
\end{equation}
See \cite{HLM}. If we choose, as we may do, $f$ and $g$ to have equal traces, then by
removing  squares in the Hilbert symbol we obtain 
\begin{equation}
AG^{(2)} \cong \left( {(\beta(f) + 4)\beta(f) , \gamma(f,g) \over kG} \right).
\end{equation}
Finally note that if $kG  = {\bf Q}(\tr(G))$ we can simplify still further to 
\begin{equation}
AG^{(2)} \cong \left( { \beta(f) , \gamma(f,g) \over kG } \right).
\end{equation}
With these representations of $AG^{(2)}$, condition 3. of Theorem 4.3 can be 
readily checked using Lemma 4.1. Thus in summary we have the following result.
\begin{theorem}
Let $G$ be a non-elementary finitely generated subgroup of $PSL(2,{\bf C})$
with a pair of non-commuting elements $f$ and $g$, where $f$ is not parabolic and 
$f$ and $g$ are both not of order 2. Then $G$ is discrete if 
\begin{enumerate}
\item ${\bf Q}({\em tr}(G))$ is a finite extension of ${\bf Q}$;
\item ${\em tr}(G)$ consists of algebraic integers;
\item $kG$ has one complex place or is totally real;
\item with $\displaystyle{AG^{(2)} \cong \biggl ({{a,b} \over kG}\biggr)}$ 
described 
at {\rm (4.4)}, {\rm (4.5)}, {\rm (4.6)}, then $\sigma(a)$ and $\sigma(b)$ are
negative for each embedding $\sigma$ of $kG$ $\neq id,{\bf c}$.
\end{enumerate}
\end{theorem}

\section {Two-generator groups}

Here we specialize the discreteness theorems in the previous sections to
the case where $G$ is a two-generator group and, in particular, where both
generators are elliptic. In these cases, both the invariant field and the 
invariant quaternion algebra will be readily described in terms of the 
parameters of the group.

As stated above, $kG = {\bf Q}({\rm tr}(G^{(2)}))$ is an invariant of the
commensurability class. In fact, it is shown in \cite{R2} that the field $kG$ 
coincides with the field 
\[ {\bf Q}( \{{\rm tr}^2(g):g \in G\})={\bf Q}(\{\beta(g):g \in G\}). \] 
Actually, it is not difficult to see that if $G$ is generated by elements $g_1,
g_2, \ldots ,g_n$ with  $\beta(g_i)\neq -4$ for $i=1,2,\ldots,n$, then
\begin{equation}
kG = {\bf Q}({\tr}(G^{sq}))
\end{equation}
where $G^{sq}=\langle g_1^2,g_2^2,\ldots,g_n^2 \rangle$;  see \cite{HLM}.
For two-generator groups this has the following consequence.

\begin{lemma}
Let  $G=\langle f,g \rangle$ be a subgroup of $SL(2,{\bf C})$ with
$\beta(f)\neq -4$ and $\beta(g)\neq -4$. Then
\begin{equation}
kG = {\bf Q}(\beta(f),\beta(g),\beta(fg^{-1})-\gamma(f,g)).
\end{equation}
\end{lemma}
{\bf Proof.} The trace of any element in a two-generator group $\langle
\phi,\psi \rangle$  group is given by a polynomial with integer coefficients in
${\tr}(\phi), {\tr}(\psi), {\tr}(\phi\psi^{-1}) $.  See, for example
\cite{Hor}.  Thus ${\bf Q}({\tr}(G^{sq}))$ is ${\bf Q}({\tr}(f^2),{\tr}(g^2),{\tr}(f^2
g^{-2}))$. Then since
\[{\tr}(f^2) = \beta(f)+2, \\\ {\tr}(g^2) = \beta(g)+2,\]
the desired conclusion follows from the identity
\[{\tr}(f^2g^{-2})=\beta(fg^{-1})+2-\gamma(f,g). \\\ \Box \]

The following corollary will be of use to us.

\begin{corollary}
Let  $G=\langle f,g \rangle$ be a group with $\beta(f) \neq -4$ and
$\gamma(f,g)\neq 0,\beta(f)$. If $G_1=\langle f,gfg^{-1} \rangle$, then
\begin{equation} 
kG_1 = {\bf Q}(\beta(f),\gamma(f,g)).
\end{equation}
\end{corollary}
{\bf Proof.} Applying Lemma 5.2 in the case where the generators are $f$ and
$gfg^{-1}$ and $f$ as above, we see that 
\begin{equation}
kG_1 = {\bf Q}(\beta(f),\beta(fgf^{-1}g^{-1})-\gamma(f,gfg^{-1})).
\end{equation}
An easy calculation shows
$$\beta(fgf^{-1}g^{-1}) - \gamma(f,gfg^{-1}) = \gamma(f,g)(\beta(f) + 4).$$
Next the hypotheses on $\beta(f)$ and $\gamma(f,g)$ imply that $f$ and 
$gfg^{-1}$ do not commute. Thus $kG_1 = {\bf Q}(\gamma(f,g),\beta(f))$. $\Box$

\begin{corollary}  Let $G=\langle f,g \rangle$ be a group with $\beta(f)
\neq -4$, $\beta(g) = -4$ and $\gamma(f,g)\neq 0,\beta(f)$. Then
\[ kG={\bf Q}(\gamma(f,g),\beta(f)).\]
\end{corollary}
{\bf Proof.}
Since $g$ has order two, $G_1=\langle f,gfg^{-1} \rangle$ has index two in
$\langle f,g \rangle$. Thus $kG_1=kG$ being an invariant of the commensurability
class and the result follows from Corollary 5.4. $\Box$

\bigskip
The hypotheses of Corollary 5.7 imply that the groups $G=\langle f,g
\rangle$ and $G_1=\langle f,gfg^{-1} \rangle$ are simultaneously
discrete or nondiscrete.  For two-generator groups $G$ we have also
identified $kG$ in terms of the parameters of $G$. We have therefore
proved the following result.

\begin{corollary}
Let $G=\langle f,g \rangle$ be a group with {\rm
par}$(G)=(\gamma,\beta,-4)$ where $\beta \neq -4,0$
and $\gamma\neq 0,\beta$. Then $kG = {\bf Q}(\gamma,\beta)$ and 
\begin{equation}
AG^{(2)} = \biggl({{\beta(\beta +4),\gamma (\gamma - \beta)}\over {\bf
Q}(\gamma,\beta)} \biggr).
\end{equation}
\end{corollary}
{\bf Proof.} Again $G_1= \langle f,gfg^{-1}\rangle$ is a subgroup of index
two in $G$ with
\[ {\rm par}(G_1)=(\gamma(\gamma-\beta),\beta,\beta).\]
Our hypotheses imply this group is nonelementary and $kG_1=kG={\bf
Q}(\gamma,\beta)$. The Hilbert symbol follows from (4.5).  $\Box$

\bigskip

From Theorem 4.7, we can now deduce as a corollary the following sufficient 
condition for a two-generator group to be a subgroup of an arithmetic Kleinian group,
and hence discrete.

\begin{theorem}
Let $G=\langle f,g \rangle$ be a subgroup of $PSL(2,{\bf C})$ with
{\rm par}$(G) =(\gamma,\beta,-4)$ where $\beta \neq 0$ and $\gamma \neq 0,\beta$. Then
$G$ is a discrete group if
\begin{enumerate}
\item  $\gamma$ and $\beta$ are algebraic integers;
\item  the field $kG={\bf Q}(\gamma,\beta)$ has exactly one complex place
or is totally real;
\item  if $kG$ has a complex place, then $-4<\sigma(\beta)<0$ and
$\sigma(\gamma(\gamma-\beta))<0$ for all real embeddings $\sigma$;
\item  if $kG$ is totally real, then $-4<\sigma(\beta)<0$ and
$\sigma(\gamma(\gamma-\beta))<0$ for all real embeddings $\sigma \neq id$.
\end{enumerate}
Indeed $G$ is a subgroup of an arithmetic Fuchsian or Kleinian group with
invariant trace field ${\bf Q}(\gamma,\beta)$ where the Hilbert symbol of the
invariant quaternion algebra is  
\[ \biggl({{(\beta+4)\beta,\gamma(\gamma-\beta)}\over {{\bf Q}(\gamma,\beta)}}
\biggr).\] 
\end{theorem}
{\bf Proof.}  The traces of $f,fg$ satisfy monic quadratic polynomials 
whose coefficients lie in ${\bf Z}[\gamma,\beta]$. Since the
traces of elements in $G$ are polynomials with integer coefficients in tr$(f)$,
tr$(g)$, tr$(fg)$, the first hypothesis implies that all traces of elements of $G$
are algebraic integers. Therefore the elements of $G^{(2)}$ also have algebraic
integer traces. By Corollary 5.8, the Hilbert symbol for $AG^{(2)}$ is as claimed,
and so by the last pair of assumptions we see that $AG^{(2)}$ is ramified at all
real places, except possibly in the case when $kG$ is totally real.

By our previous discussion it follows that ${\cal O}G^{(2)}$ is an order of
$AG^{(2)}$, and therefore ${{\cal O}G^{(2)}}^1$  yields an arithmetic Kleinian or
Fuchsian group via the construction discussed previously.  Thus $G^{(2)}$ is
discrete and hence so is $G$.

In the proof we constructed the arithmetic group ${{\cal O}G^{(2)}}^1$
which contains $G^{(2)}$. To get an arithmetic group containing $G$ consider the
group 
\[{\rm Norm}({\cal O}G^{(2)}) = \{x \in AG^{(2)} : x{\cal O}G^{(2)}x^{-1} =
{\cal O}G^{(2)}\}.\] 
The image of this group is an arithmetic subgroup of $\PSL(2,{\bf C})$
(see \cite{Bo} for example). As $G^{(2)}$ is normal in $G$, it
follows that $G$ is a subgroup of an arithmetic Kleinian group. $\Box$
\bigskip

In the cases where $f$ is primitive elliptic of order $n$ and $g$ is elliptic 
of order 2, then $\beta = -4 \sin^2 \pi/n$ is a totally real algebraic 
integer so $\gamma$ must satisfy a monic polynomial in ${\bf Z}[\beta][z]$
to be an algebraic integer. For all embeddings $\sigma$ of ${\bf Q}(\gamma,
\beta)$, $-4 < \sigma(\beta) \leq \beta < 0$ and so Theorem 2.6 is an 
immediate corollary of Theorem 5.10.

\bigskip

We now refine these results to obtain more easily applicable criteria for
discreteness. We first recall some basic results concerning field extensions.

Suppose that ${\bf K}$ is a finite extension of ${\bf Q}$, that $\gamma$ is
algebraic over ${\bf K}$ and that $p(z)$ is the minimum polynomial of $\gamma$
over ${\bf K}$ with deg$(p) = m > 1$.  Let $\sigma:{\bf K} \rightarrow {\bf C}$ be
an embedding.  Then there are exactly $m$ embeddings $\tau:{\bf K}(\gamma)
\rightarrow {\bf C}$ such that $\tau|{\bf K}=\sigma$ and these are uniquely
determined by $\tau(\gamma) = \gamma'$ where $\gamma'$ runs through the
roots of $\sigma(p(z))$. See for example \cite{Ga}.

\begin{lemma}
If ${\bf K}(\gamma)$ has exactly one complex place, then ${\bf K}$ must be
totally real, i.e. every embedding of ${\bf K}$ into ${\bf C}$ must be real.
\end{lemma}
{\bf Proof.}
If $\sigma:{\bf K} \rightarrow {\bf C}$ is a complex embedding,  so is
$\bar \sigma$ and $\bar \sigma \neq \sigma$.  Then by the above remarks, there
exist 2$m$ embeddings of ${\bf K}$ which are complex. $\Box$

\begin{theorem}
Let $\gamma$ and $\beta$ be algebraic integers such that
$\gamma \not \in {\bf Q}(\beta) {\subset} {\bf Q}(\gamma,\beta)$ and let
$p(z,\beta)$ be the minimum polynomial of $\gamma$ over ${\bf Q}(\beta)$.
Then ${\bf Q}(\gamma,\beta)$ has one complex place if and only if 
\begin{enumerate}
\item $\beta$ is totally real, i.e. all of the Galois conjugates of
$\beta$, $\beta$=$\beta_1$, $\beta_2$, \ldots, $\beta_n$ are real;
\item the polynomial $q(z)=p(z,\beta_1)p(z,\beta_2)\cdots p(z,\beta_n)$
has exactly one pair of complex conjugate roots and all other roots real.
\end{enumerate}
\end{theorem}
{\bf Proof.}
For $j=1,2,\ldots,n$ let $\sigma_j$ denote the embedding of ${\bf
Q}(\beta)$ into ${\bf C}$ defined
by $\sigma_j(\beta) = \beta_j$ and suppose that $\tau$ is an embedding of
${\bf Q}(\gamma,\beta)$
into ${\bf C}$.  Then $\tau$ is determined by the images $\tau(\gamma),
\tau(\beta)$.  Now
$\tau(\beta) = \sigma_j(\beta)$ for some $j$ in which case $\tau(\gamma)$
is a root of
$\sigma_j(p(z,\beta)) = p(z,\beta_j)$.  Thus conditions 1. and 2. imply that there
are precisely two complex embeddings ${\bf Q}(\gamma,\beta)$ into ${\bf C}$.

Conversely if ${\bf Q}(\gamma,\beta)$ has exactly one pair of complex
embeddings into ${\bf C}$, then by Lemma 5.11 each embedding of ${\bf Q}(\beta)$
into ${\bf C}$ must be real and so an embedding $\tau$ of ${\bf Q}(\gamma,\beta)$
into ${\bf C}$ will be complex precisely when a root of $\sigma_j(p(z,\beta)) =
p(z,\beta_j)$ is complex. $\Box$

\bigskip

Note that with $\gamma,\beta$ as in Theorem 5.12, ${\bf Q}(\gamma)$
must have at least one complex place and so by Lemma 5.11, ${\bf
Q}(\gamma) = {\bf Q}(\gamma,\beta)$.

Now specialize again to the cases where $G = \langle f,g \rangle$ with 
$f$ a primitive elliptic element of order $n \geq 3$, $g$ elliptic of order 2
and $\gamma \neq 0,\beta$. Thus ${\bf Q}(\beta)$ is totally real and the Galois 
conjugates $\beta_k = -4\sin^2(k\pi/n)$ where $(k,n) = 1$ and $1 \leq k 
\leq n/2$ so that $\beta_k = \beta(f^k)$. Note that $-4 < \beta_k < \beta 
< 0$ for $k \neq 1$. Thus from Theorem 5.10 and 5.12, we obtain the following result.

\begin{theorem}
Let $G=\langle f,g \rangle$ be a subgroup of $PSL(2,{\bf C})$ with $f$ a primitive
elliptic of order $n \geq 3$, $g$ elliptic of order $2$ and $\gamma(f,g) \neq
0,\beta(f)$. Then $G$ is a discrete subgroup of an arithmetic group if
\begin{enumerate}
\item $\gamma(f,g)$ is the root of a polynomial $p(z,\beta(f))$ in $z$ and
$\beta(f)$ with integer coefficients;
\item the roots of $p(z,\beta(f^k))$ are real and lie in the interval 
$(\beta(f^k),0)$ for all $k$ such 
that $(k,n) = 1$ and $2 \leq k \leq n/2$;
\item if $\gamma(f,g)$ and $\bar \gamma(f,g)$ are complex, then all other
roots of $p(z,\beta(f))$ lie in the interval $(\beta(f),0)$;
\item if $\gamma(f,g)$ is real, then all other roots of $p(z,\beta(f))$ lie
in the interval $(\beta(f),0)$.   
\end{enumerate}
\end{theorem}

Finally if $n = 3,4$ or $6$, then $\beta \in {\bf Z}$ and we have the following result. 

\begin{theorem}
Let $G=\langle f,g \rangle$ be a subgroup of $PSL(2,{\bf C})$ with $f$ a
primitive elliptic of order $n=3,4$ or $6$, $g$ an elliptic of order 2 and
$\gamma(f,g) \neq 0,\beta(f)$.  Then  $G$ is a discrete subgroup of an arithmetic
group if  
\begin{enumerate}
\item $\gamma(f,g)$ is the root of a polynomial $p(z)$ with integer
coefficients;
\item if $\gamma(f,g)$ and $\bar \gamma(f,g)$ are complex, then all other
roots of $p(z)$ lie in the
interval  $(\beta(f),0)$;
\item if $\gamma(f,g)$ is real, then all other roots of $p(z)$ lie in the
interval $(\beta(f),0)$.
\end{enumerate}
\end{theorem}

\section{Examples}

In this section we use the discreteness criteria in \S 5 to identify a number of
interesting two-generator subgroups of arithmetic groups. Some of these were shown in
\cite{GM3}, \cite{GM4} and \cite{GM5} to be extremal for geometric properties such as
minimum axial distance, minimum co-volume and the Margulis constant. A number of these
groups are of finite co-volume themselves and are therefore arithmetic; recall that
our methods identify the invariant trace field and quaternion algebra.

The groups we consider are generated by a pair of elliptic elements $f$ and
$g$ where $f$ is of order $n \geq 3$  and $g$ of order $2$. We present these in
tabular form for the cases $n=3,4,5,6,7$. By virtue of Theorems 5.13 and 5.14 all of
these are subgroups of  arithmetic groups, and hence discrete. We discuss in \S 7
and \S 8 their relations with small volume hyperbolic 3-orbifolds and in \S 10 we
determine the groups in which the generator $f$ is a simple elliptic.

We begin with a pair of preliminary remarks. Suppose that $(\gamma,\beta,-4)$ are
the parameters of a discrete two--generator group $G=\langle f,g \rangle$ with $f$
of order $n$. Then the Lie product $G^*$ of $G$ \cite{J} and the conjugations $\bar
G$ and $\bar G^*$ of $G$ and $G^*$ by reflection in the real axis are discrete
two-generator groups with 
\[ {\rm par}(G^*)=(\beta - \gamma,\beta,-4), 
\\\ {\rm par}(\bar G)=(\bar \gamma,\beta,-4), 
\\\ {\rm par}(\bar G^*)=(\beta - \bar \gamma,\beta,-4).\] 
See, for example, \cite{GM3}.  Moreover, if $G$ is a subgroup of an arithmetic
group, then the same is true of $G^*$, $\bar G$ and $\bar G^*$.  Hence it is
sufficient to consider only those values of
$\gamma$ for which
\[ -2\sin^2(\pi/n) = \beta/2 \leq {\rm Re}(\gamma),\\\ \ 0 \leq {\rm Im}(\gamma). \]

Next if $\gamma \neq 0,\beta$, then by Lemma 2.4 the hyperbolic distance $\delta
=\delta(f,g)$ between the axes of the generators $f$ and $g$ of $G$ is given by
\begin{equation}
\cosh(2 \delta)=\frac{|\gamma - \beta|+|\gamma|}{|\beta|}.
\end{equation}	
In particular, $\delta = 0$ if and only if $\langle f,g \rangle$ is elementary.
Hence in our case, $\delta$ is a measure of how much $\langle f,g \rangle$
differs from an elementary group. We see from (6.1) that the corresponding axial
distances for $G^*$, $\bar G$ and $\bar G^*$ are equal to that for $G$.

We now give for $n=3,4,5,6,7$ tables of groups $G_{n,i}$ generated by
elliptics of orders $n$ and $2$ together with
\begin{enumerate}
\item the approximate values of their commutator parameters $\gamma_{n,i}$,
\item the minimum polynomial $p_{n,i}$ for $\gamma_{n,i}$ over ${\bf
Q}(\beta)$,
\item the approximate distance $\delta_{n,i}$ between the axes of the
generators of $G_{n,i}$.
\end{enumerate}
When $n=3,4,6$, $p_{n,i}$ coincides with the minimum polynomial $q_{n,i}$ for
$\gamma_{n,i}$ over the field ${\bf Q}$ considered in the next section.

These groups were identified by the disk covering procedure described in \S 2.  The
resulting diagrams, showing the only possible values for the parameter $\gamma_{n,i}$ 
for the cases where $n=3,4,5,6$, are given at the end of this paper in \S 11. For an
account of the calculations on which these diagrams are based, see \cite{GM6}. 

We begin with tables for the groups $G_{3,i}$ and $G_{4,i}$.

\bigskip

\begin{center}
\begin{tabular}{|c|c|c|c|}
\multicolumn{4}{c}{\bf Table 1 - Groups {\boldmath $G_{3,i}$}}\\
\hline
$i$&$\gamma_{3,i}$ &$p_{3,i}$ & $\delta_{3,i}$\\
\hline
$ 1 $&$ -1 $&$ z+1 $&$ 0 $\\
\hline
$ 2 $&$ -.3819 $&$ z^2+3z+1  $&$ 0 $\\
\hline
$ 3 $&$ -1.5+.6066i $&$ z^4+6z^3+12z^2+9z+1 $&$ .1970 $\\
\hline
$ 4 $&$ -.2118+.4013i $&$ z^4+5z^3+7z^2+3z+1 $&$ .2108 $\\
\hline
$ 5 $&$ -.5803+.6062i $&$ z^3+4z^2+4z+2 $&$ .2337 $\\
\hline
$ 6 $&$ -1.1225+.7448i $&$ z^3+5z^2+8z+5 $&$ .2448 $\\
\hline
$ 7 $&$ -.3376+.5622i $&$ z^3+3z^2+2z+1 $&$ .2480 $\\
\hline
$ 8 $&$ -.9236+.8147i $&$ z^4+5z^3+8z^2+6z+1 $&$ .2740 $\\
\hline
$ 9 $&$ -1.5+.8660i $&$ z^2+3z+3 $&$ .2746 $\\
\hline
$ 10 $&$ -.7672+.7925i $&$ z^3+4z^2+5z+3 $&$ .2770 $\\
\hline
$ 11 $&$ .0611+.3882i $&$ z^4+5z^3+6z^2+1$&$ .2788 $\\
\hline
$ 12 $&$ .2469 $&$ z^3+4z^2+3z-1 $&$ .2831 $\\
\hline
$ 13 $&$ -.5284+.7812i $&$ z^6+8z^5+24z^4+35z^3+28z^2+12z+1 $&$ .2944 $\\
\hline
$ 14 $&$ -.1153+.5897i $&$ z^3+3z^2+z+1 $&$ .2970 $\\
\hline
\end{tabular}
\bigskip

\begin{tabular}{|c|c|c|c|}
\multicolumn{4}{c}{\bf Table 2 - Groups {\boldmath $G_{4,i}$}}\\
\hline
$i$&$\gamma_{4,i}$ &$p_{4,i}$ & $\delta_{4,i}$\\
\hline
$ 1 $&$ -1 $&$ z+1 $&$ 0 $\\
\hline
$ 2 $&$ -.5+.8660i $&$ z^2+z+1  $&$ .4157 $\\
\hline
$ 3 $&$ -.1225+.7448i $&$ z^3+2z^2+z+1 $&$ .4269 $\\
\hline
$ 4 $&$ -1+i $&$ z^2+2z+2 $&$ .4406 $\\
\hline
$ 5 $&$ -.6588+1.1615i $&$ z^3+3z^2+4z+3 $&$ .5049 $\\
\hline
$ 6 $&$ .2327+.7925i $&$ z^3+z^2+1 $&$ .5225 $\\
\hline
$ 7 $&$ -.2281+1.1151i $&$ z^3+2z^2+2z+2 $&$ .5297 $\\
\hline
$ 8 $&$ .4196+.6062i $&$ z^3+z^2-z+1 $&$ .5297 $\\
\hline
$ 9 $&$ i $&$ z^2+1 $&$ .5306 $\\
\hline
$ 10 $&$ .6180 $&$ z^2+z-1 $&$ .5306 $\\
\hline
$ 11 $&$ -1+1.2720i $&$ z^4+4z^3+7z^2+6z+1$&$ .5306 $\\
\hline
$ 12 $&$ -.4063+1.1961i $&$ z^4+3z^3+4z^2+4z+1 $&$ .5345 $\\
\hline
$ 13 $&$ .7881+.4013i $&$ z^4+z^3-2z^2+1 $&$ .6130 $\\
\hline
\end{tabular}
\end{center}
\bigskip

Tables 3 and 5 contain the groups for $n=5$ and $7$.  Here $p_{n,i}$ is a
polynomial in $z$ and $\beta$ to which Theorem 5.13 applies.

\bigskip
\begin{center}
\begin{tabular}{|c|c|c|c|}
\multicolumn{4}{c}{\bf Table 3 - Groups {\boldmath $G_{5,i}$}}\\
\hline
$i$&$\gamma_{5,i}$ &$p_{5,i}$ & $\delta_{5,i}$\\
\hline
$ 1 $&$ -.3819 $&$ z-\beta-1 $&$ 0 $\\
\hline
$ 2 $&$-.6909+.7228i $&$ z^2-\beta z+1  $&$ .4568 $\\
\hline
$ 3 $&$ .1180+.6066i $&$ z(z-\beta-1)^2-\beta-1 $&$ .5306 $\\
\hline
$ 4 $&$ -.1909+.9815i $&$ z^2-(\beta+1)z+1 $&$ .6097 $\\
\hline
$ 5 $&$ .6180 $&$ z-\beta-2 $&$ .6268 $\\
\hline
$ 6 $&$ .2527+.8507i $&$ z(z-\beta-1)^2+1 $&$ .6514 $\\
\hline
$ 7 $&$ -.6909+1.2339i $&$ z^2-\beta z+2 $&$ .6717 $\\
\hline
$ 8 $&$ -.3819+1.2720i $&$ z^3-(2\beta +1)z^2+(\beta^2+\beta+2)z-2\beta-1
$&$ .6949 $\\
\hline
$ 9 $&$ .1180+1.1696i $&$ z^3-(2\beta+2)z^2+(\beta^2+2\beta+2)z-\beta $&$
.7195 $\\
\hline
$ 10 $&$ -.0817+1.2880i $&$ z^4-(2\beta+1)z^3+(\beta^2+\beta+2)z^2-2\beta
z+1 $&$ .7273 $\\
\hline
$ 11 $&$ .6180+.7861i $&$ z^3-(2\beta+3)z^2+(\beta^2+3\beta+2)z+1 $&$ .7323 $\\
\hline
$ 12 $&$ .8776+.5825i $&$ z(z-\beta)(z-\beta-2)^2+1 $&$ .7725 $\\
\hline
\end{tabular}
\bigskip

\begin{tabular}{|c|c|c|c|}
\multicolumn{4}{c}{\bf Table 4 - Groups {\boldmath $G_{6,i}$}}\\
\hline
$i$&$\gamma_{6,i}$ &$p_{6,i}$ & $\delta_{6,i}$\\
\hline
$ 1 $&$ -.5+.8660i $&$ z^2+z+1 $&$ .6584 $\\
\hline
$ 2 $&$ i $&$ z^2+1  $&$ .7642 $\\
\hline
$ 3 $&$  .5+.8660i $&$ z^2-z+1 $&$ .8314 $\\
\hline
$ 4 $&$ -.5+1.3228i $&$ z^2+z+2 $&$ .8500 $\\
\hline
$ 5 $&$ -.2150+1.3071i $&$ z^3+z^2+2z+1 $&$ .8539 $\\
\hline
$ 6 $&$ .3411+1.1615i $&$ z^3+z+1 $&$ .8786 $\\
\hline
$ 7 $&$ 1 $&$ z-1 $&$ .8813 $\\
\hline
$ 8 $&$ .8774+.7448i $&$ z^3-z^2+1 $&$ .9106 $\\
\hline
\end{tabular}
\bigskip

\begin{tabular}{|c|c|c|c|}
\multicolumn{4}{c}{\bf Table 5 - Groups {\boldmath $G_{7,i}$}}\\
\hline
$i$&$\gamma_{7,i}$ &$p_{7,i}$ & $\delta_{7,i}$\\
\hline
$ 1 $&$ .2469 $&$ z-\beta-1 $&$ .5452$\\
\hline
$ 2 $&$ -.3765+.9264i $&$ z^2-\beta z+1  $&$ .8162 $\\
\hline
$ 3 $&$  1.2469 $&$ z-\beta-2 $&$ 1.0704 $\\
\hline
\end{tabular}
\end{center}
\bigskip

We conclude this section with two examples to illustrate how the above tables are
constructed.  The groups we consider,  $G_{3,2}$ and $G_{5,2}$, are those where the
minimum distances between the axes of elliptics of orders 3 and 2 and orders 5 and 2 are
realized.

In the first example, $\beta = \beta(f) = -3$ and the disk covering procedure shows
that the polynomial
 \[z(z+3)(z^2+3z+3)(z^4+6z^3+12z^2+9z+1)^2(z^4+6z^3+12z^2+9z+3), \]
corresponding to the word $hfh^{-1}fhfh^{-1}$ where
$h=gfg^{-1}fgf^{-1}g^{-1}$, maps
\[\gamma_{3,2}=-3/2+i\sqrt{(2\sqrt{5}-3)/4} \approx -1.5+.606658i\]
onto $0$. Then \[p_{3,2}(z) = z^4+6z^3+12z^2+9z+1\] is the minimum
polynomial for  $\gamma_{3,2}$
and it is easy to check that $p_{3,2}$ has real roots $-.13324$, $-2.86676$
which lie in the interval
$(-3,0)$.  Hence we can apply Theorem 5.14 to conclude that $G_{3,2}$ is a
subgroup of an arithmetic
group and hence discrete.

In the second case,  
\[\beta = \beta(f) = (\sqrt{5}-5)/2, \;\;\; \beta(f^2)=(-\sqrt{5}-5)/2\] 
and from the disk covering argument we find that  \[p_{5,2}(z,\beta) = z^2-\beta
z-\beta\] is the minimum polynomial for
\[\gamma_{5,2} = (\sqrt{5}-5)/4+i \sqrt{(5 \sqrt{5}-7)/8} \approx
-.690983+.722871i \]
over ${\bf Q}(\beta)$.  Then  $\gamma_{5,2}$ and $\bar \gamma_{5,2}$ are the
roots of $p_{5,2}(z,\beta(f))$, the roots  $-.301522$ and $-3.31651$ of
$p_{5,2}(z,\beta(f^2))$ are both real and the desired conclusion follows from
Theorem 5.13.

\section{ Small volume arithmetic orbifolds}

We now discuss the minimal co-volume arithmetic Kleinian groups that the groups
$G_{n,i}$ embed in. We deal with the case $n = 3$ in some detail and in \S 8 merely
state the salient information for the others in tabulated form.

We begin with some comments on the more detailed arithmetic structure associated to
these groups. For this we need some additional results and terminology,
for which we refer the reader to 
\cite{V}.

Let $B$ be a quaternion algebra over a number field $k$. The isomorphism 
class of $B$ is determined by those places of $k$, both Archimedean and non-Archimedean, 
at which $B$ is ramified. The number of ramified places is always even. 
The non-Archimedean places correspond to the prime ideals in $R$, the 
ring of integers in $k$. If ${\cal O}$ is an order in $B$, then its 
{\em discriminant} $d({\cal O})$ is an ideal of $R$. Furthermore, 
${\cal O}$ will be a maximal ideal precisely when $d({\cal O})$ is the 
product of those ideals at which $B$ is ramified. If ${\cal M}$ is an 
order of $B$ such that ${\cal M} \subset {\cal O}$, then 
$d({\cal O}) \mid d({\cal M})$. In the cases where ${\cal O}$ is a free 
$R$-module with basis $\{f_1,f_2,f_3,f_4\}$, the discriminant can be 
determined from the formula 
$$ d({\cal O})^2 = \langle {\rm Det}({\tr_B}(f_if_j)) \rangle$$
where $\tr_B : B \rightarrow k$ is the reduced trace.

\medskip

The following lemma will prove useful in our calculations.

\begin{lemma}
Let $G$ be a finite co-volume Kleinian group such that ${\rm tr}(G)$ consists of
algebraic integers and let $R$ denote the ring of integers in ${\bf Q}({\rm tr}(G))$.
If $\langle f,g \rangle$ is a nonelementary subgroup of $G$, then $R[1,f,g,fg]$ is
an order of $AG$. 
\end{lemma}
{\bf Proof.}
Recall that $1,f,g,fg$ span $AG$ over the trace field. Thus 
since $R[1,f,g,fg]$ contains a basis of $AG$, is finitely generated and
contains $R$, it suffices to prove that all products of the basis elements can be
expressed as $R$-combinations of the basis elements. There are several obvious ones
and these together with the following identities prove the lemma.
\begin{eqnarray}
f^2 & = & {\tr} (f)f - 1, \nonumber \\
f^{2} g & = & {\tr} (f)fg - g,\nonumber \\
(fg)^2 & = & {\tr} (fg)fg - 1, \nonumber \\
fgf & = &-{\tr} (g)1+{\tr} (fg)f+g, \nonumber \\
gf+fg & = &({\tr} (fg)-{\tr}(f){\tr}(g))1+{\tr}(g)f+{\tr}(f)g. \;\;\; \Box \nonumber
\end{eqnarray}

We adopt some notation that will be used in the next two sections.
For $n=3$,$4$,$5$,$6$,$7$ let $kG_{n,i}$ and $A_{n,i}$ denote the invariant
trace field and invariant quaternion algebra of $G_{n,i}$.

We now fix attention on the case of $n=3$ and, in what follows, suppress the
subscript $3$ for convenience of notation. Next, referring to Table 1 in the cases 
$i$=1, 2 and 12, since the corresponding commutator parameter is real, one checks
easily that these groups are the (2,3,4) spherical triangle group, the (2,3,5)
spherical triangle group and the (2,3,7) Fuchsian triangle group.  Also in the case
$i$=13, the invariant trace field has degree 6 and we will not make any further
comment on this case.  Thus we concentrate on the cases when the commutator parameter
is not real and the relevant polynomials have degree 2, 3 and 4.

Let ${\cal O}_i$ denote the  $R_{kG_i}$-submodule of $A_i$ generated over
$R_{kG_i}$ by the elements 
\[\{1, g_if_ig_i^{-1}, f_i^{-1}, g_if_ig_i^{-1}f_i^{-1}\};\] 
note that since $f_i$ has order $3$, $f_i$ and $g_if_ig_i^{-1}$ lie in $G_i^{(2)}$.
Lemma 7.1 shows that ${\cal O}_i$ is an order of $A_i$.

By Theorem 5.10, the algebras $A_i$ have the following Hilbert symbols,
\[ \biggl({{-3,(\gamma_i+3)\gamma_i}\over kG_i} \biggr). \]
It is an easy calculation using the comments above to compute
the discriminant of ${\cal O}_i$, namely
\[ d({\cal O}_i) = \langle{\tr}([f_i^{-1},g_if_ig_i^{-1}])-2 \rangle =
\langle \gamma_i(\gamma_i+3) \rangle. \]

We first deal with the case of the quartics and cubics.

\begin{lemma}
For each $i \neq 10$,  corresponding to a cubic or quartic,
${\cal O}_i$ is a maximal order.
\end{lemma}
{\bf Proof.} Using the earlier remarks we need to compute the ramification
of these quaternion
algebras and compare them with $\gamma_i(\gamma_i+3)$.
Since each ramified prime must divide $d({\cal O}_i)$, it is expeditious to 
calculate $d({\cal O}_i)$ in the first instance.

In all the quartic cases, a simple calculation shows that $\gamma_i(\gamma_i 
+ 3)$ is a unit, so that $d({\cal O}_i) = R_{kG_i}$, from which it 
follows that $A_i$ cannot have any finite ramification and ${\cal O}_i$ 
is maximal.

A variation of this argument works in the cubic cases, $i \neq 10$. For 
the norm of $\gamma_i(\gamma_i + 3)$ in each of these cases is a rational prime
so that $d({\cal O}_i)$ is a prime ideal ${\cal P}_{p_i}$. As the cardinality 
of the set of places ramified in any $A_i$ is even and as each $A_i$ here 
is ramified at exactly one real place, $A_i$ must be ramified at at least 
one prime ideal. Thus it is ramified at precisely the one ideal 
${\cal P}_{p_i}$ and ${\cal O}_i$ is maximal. For later computations, we 
note that $p_5 = 2,p_6 = 5,p_7 = 5,p_{14}= 2$. $\Box$

\medskip

For the case $i=10$ we argue as follows.  The discriminant of $kG_{10}$ is -31 and by
the fact that there is a unique cubic with one complex place and discriminant -31, an
alternative description of $kG_{10}$ is $kG_{10} = {\bf Q}(u)$ where $u^3 + u + 1 = 0$;
cf \cite{PZW}.  The ring of integers in this field has a unique prime ideal ${\cal P}_3$
of  norm 3 and ${\cal P}_3 = < \gamma_{10}>$. Now $\gamma_{10} + 3$ also has 
norm 3 so that $\gamma_{10}( \gamma_{10} + 3) = v \gamma_{10}^2$ 
where $v$ is a unit. Thus $d({\cal O}) = {\cal P}_3^2$. As in the above proof
$A_{10}$ must be ramified at at least one prime and so must be 
ramified at precisely the one finite place corresponding to ${\cal P}_3$.
Note that ${\cal O}_{10}$ is not a maximal order in this case. $\Box$

\medskip

We now compute maximal discrete groups in which the groups $G_i$, with $i$ as
indicated above, are subgroups. To do this we make use of the description of maximal
groups due to Borel \cite{Bo}.  First we consider only cubics and
quartics and we assume that $i \neq 10$. The case $i=10$ will be dealt with
separately.

Borel's classification of maximal groups in the commensurablity class
says that if $\cal O$ is a maximal order of the quaternion algebra
$B$, then the group $P\rho({\rm Norm}({\cal O}))$ (recall the proof of
Theorem 5.10) is a maximal group of minimal co-volume in that commensurability 
class. We
denote this maximal group by $\Gamma_{\emptyset,\emptyset}$ in the
notation of \cite{Bo}.
Borel actually shows there are infinitely many maximal groups in the
commensurability class of an arithmetic group; see the case of $i = 10$ below. 

As in the proof of Theorem 5.10, since $g_i$ normalizes $\langle
f_i,g_if_ig_i^{-1}\rangle$, we deduce that $G_i < P\rho(\Norm({\cal O}_i))$.  By Lemma
7.2, for $i \neq 10$ ${\cal O}_i$ is maximal and we obtain the following result.

\begin{corollary}
For $i \neq 10$, $G_i$ is a subgroup of the group
 $P\rho({\rm Norm}({\cal O}_i)) =  \Gamma_{\emptyset,\emptyset}^{(i)}$.
\end{corollary}

For the case $i = 10$, we utilize more of Borel's construction of maximal
groups referred to above. We briefly recall the relevant points.  We recommend the
reader have \cite{Bo} at hand. In what follows $B$ is a quaternion algebra over the
number field $k$ with one complex place, which is ramified at all real places. Let $V$
denote the set of all finite places of $k$, and  ${\rm Ram}_f(B)$ denote the finite
places which ramify $B$.

Borel proves that to any pair $S$ and $S'$ of finite (possibly empty),
disjoint subsets
of $V \setminus {\rm Ram}_f(B)$, one can associate a group $\Gamma_{S,S'}$.
Moreover
it is shown in \cite{Bo} any arithmetic Kleinian group $\Gamma$ in the
commensurability class determined
by $B$ is conjugate to a subgroup of some $\Gamma_{S,S'}$. We will not make
use of this explicitly here.

Let $\nu \in V \setminus {\rm Ram}_f(B)$ and denote by
${\cal C}_{\nu}$ be the Bruhat-Tits
Building of ${\rm SL}(2,k_{\nu})$, which in these cases is a tree. The vertices of 
this tree are maximal orders of $M(2,k_{\nu})$, and the ${\rm SL}(2,k_{\nu})$ action on
the vertices forms 2 orbits. The groups $\Gamma_{S,S'}$ are described in terms
of vertex and edge stabilizers of these trees. Since $B$ is unramified at $\nu$
we get $B_{\nu} = B \otimes_k k_{\nu} \cong M(2,k_{\nu})$. Thus for any group in the
commensurability class determined by $B$ we have an induced
action on ${\cal C}_{\nu}$.

In our case, since $d({\cal O}_{10})$ is divisible only by ${\cal P}_3$,
if $\nu \in V \setminus {\cal P}_3$, then ${\cal O}_{10}$, localized at $\nu$
is maximal; see \cite{V}. Given any element of ${\rm Norm}({\cal O}_{10})$, the image of
this in $A_{10,\nu}$ stabilizes the vertex of ${\cal C}_{\nu}$ associated to ${\cal
O}_{10}$ localized at $\nu$. It now follows from this and Borel's description of maximal
groups that $P\rho({\rm Norm}({\cal O}_{10})) \subset \Gamma^{10}_{\emptyset,\emptyset}$; see
section 4 of \cite{Bo}.  Hence we can conclude that Corollary 7.3 also holds for $i =
10$.

\medskip
To compute the minimal co-volumes in the cubic and quartic cases we use the
formula of \cite{Bo}.

\medskip
\centerline{\bf Quartic cases}

\medskip
We have seen from the proof of Lemma 7.2 that in the case of $i = 3$, $4$, $8$
and $11$ the quaternion algebras $A_i$ are unramified at all finite places. It is
an easy calculation to show that the class number of these quartic fields is $1$
and further that the structure of the group of units of these fields implies that
the co-volume of the groups $\Gamma_{\emptyset,\emptyset}^{(i)}$ is
\[{|d_{kG_i}|^{3/2}\zeta_{kG_i}(2)}\over {2^7 \pi^6},\] 
where $d_{kG_i}$ is the discriminant and $\zeta_{kG_i}(s)$ is the Dedekind zeta function
of $kG_i$; cf. \cite{Bo}. Reasonably good estimates for $\zeta_{kG_i}(2)$ are obtained
by considering primes of small norm in $kG_i$ in the first few terms of the Euler
product expansion of $\zeta_{kG_i}(2)$.
                            
\begin{enumerate}
\item The polynomial $z^4 + 6z^3 + 12z^2 +9z + 1$ has discriminant $-275$, and this is
the smallest discriminant of a field with one complex place; cf. \cite{G}. As shown
in \cite{CF} for example this is the invariant quaternion algebra of the smallest volume
orientable arithmetic hyperbolic 3-orbifold. This minimal volume is approximately
$0.03905\ldots$.

\item The polynomial $z^4 +  5z^3 + 7z^2 + 3z +1$ has discriminant $-283$ and
hence the discriminant of $kG_4$ coincides with this. There is a unique quartic
field of this discriminant; cf. \cite{G}. The minimal co-volume in the commensurability
class determined by $A_4$ is approximately $0.0408$, the second smallest volume
known. This orbifold is commensurable with $(5,1)$-Dehn surgery on the figure eight knot
complement \cite{C}.

\item The other two quartics have discriminants $-491$ and $-563$, respectively. By
comparing with \cite{G}, these are the discriminants of the invariant trace-fields. The
smallest co-volume in the commensurability classes are approximately $0.1028 \ldots$ and
$0.1274\ldots$, respectively.

\end{enumerate}

\medskip
\centerline{\bf Cubic cases}

The corresponding volume formula in this case is
\[{|d_{kG_i}|^{3/2}\zeta_{kG_i}(2) (N{\cal P}-1)}\over {2^6 \pi^4},\]
using the same terminology as above. The additional term $N{\cal P}$ denotes the
norm of the prime ideal $\cal P$ ramifying the algebra in each case.

\begin{enumerate}
\item The polynomial $z^3 + 5z^2 + 8z + 5$ has discriminant $-23$ which is necessarily
the discriminant of the invariant trace-field. The minimal co-volume of a group in this
commensurability class is approximately $0.07859$. This orbifold is commensurable with
the Weeks' manifold of volume $0.942707\ldots$

There is a unique field of discriminant $-23$. Hence a non-real root of the polynomial
$z^3 + 3z^2 + 2z +1$ also generates $kG_6$.  In fact, up to conjugation, $G_7$ is a
subgroup of the same maximal group as contains $G_6$.

\item The polynomial $ z^3 +4z^2 + 5z + 3$ has discriminant $-31$ which also must
coincide with the discriminant of $kG_{10}$. Following  the proof of Lemma 7.2, we
observed  that $A_{10}$ is ramified at the place $\nu_3$ corresponding to
${\cal P}_3$.
From our comments above $G_{10}$ is a subgroup of the minimal co-volume group which
in this case is approximately $0.06596$.

(In unpublished calculations we have shown that this orbifold is commensurable with
$(8,1)$-Dehn surgery on the figure-eight knot complement, and contains the orbifold
group of $(3,0)$,$(3,0)$-Dehn surgery on the Whitehead link 
as a subgroup of index $8$. See also \cite{HLM})

\item The polynomial $z^3 + 4z^2 + 4z + 2$ has discriminant $-44$, which coincides with
the discriminant of $kG_{5}$. The smallest co-volume of a group in the commensurability
class is approximately $0.066194$.

\item The polynomial $z^3 + 3z^2 + z + 1$ has discriminant $-76$ which is the
discriminant of $kG_{14}$. The smallest co-volume in this commensurability class is
$0.1642$.

\end{enumerate}

\medskip
This leaves the case of $i = 9$ for which the relevant polynomial is quadratic.
A root of this polynomial generates the the field ${\bf Q}(\sqrt{-3})$. Moreover notice
that directly from the polynomial we can read off that $z^2 + 3z = -3$. Hence the
algebra is 
\[\biggl({{-3,-3}\over {{\bf Q}(\sqrt{-3})}} \biggr).\] 
By \cite{V}, this is isomorphic to $M(2,{\bf Q}(\sqrt{-3}))$. Indeed in this case, the
trace field of $G_{9}$ is ${\bf Q}(\sqrt{-3})$ and therefore one easily deduces from
\cite{NR} for example that $G_{9}$ is conjugate to a subgroup of $PSL(2,O_3)$ and hence
of the minimal co-volume group $\PGL(2,O_3)$. This group  has co-volume $0.08457\ldots$
and is the smallest volume of an orientable cusped hyperbolic 3-orbifold \cite{M}.

\section{Axial distance and small volumes}

In this section, we summarize calculations similar to those carried out in \S 7 for the
cases where $n= 4,5,6,7$, omitting those where the degree of the minimum polynomial
exceeds four. There are some additional complications which we comment on below, but we
make no attempt at detailing all the calculations.

In each case, the Hilbert symbol of the invariant quaternion algebra is easily computed
in terms of the parameter $\gamma_{n,i}$ as  
\begin{equation}
\biggl(\frac{-4\sin^2(2\pi/n),\gamma_{n,i}(\gamma_{n,i}+4\sin^2(\pi/n))}
{kG_{n,i}}\biggr).
\end{equation}

As in the case $n = 3$, we construct from the two-generator group $\langle f,g \rangle$
with $f^n = g^2 = 1$ a suitable order which  enables us to determine the finite
ramification of the quaternion  algebra and a small co-volume group, containing the
normalizer of the order, in which the group $G_{n,i}$ lies. 

Now for $n=3,4,5,6,7$ we give tables for the groups $G_{n,i}$ of \S 6 with
\begin{enumerate}
\item the minimum polynomial $q_{n,i}$ over ${\bf Q}$ for the commutator
parameter $\gamma_{n,i}$ of $G_{n,i}$,
\item the discriminant $d_{G_{n,i}}$ of the trace-field $kG_{n,i}$,
\item the finite places ${\rm Ram}_{f}$ where the invariant quaternion algebra
ramifies,
\item the approximate distance $\delta_{n,i}$ between the axes of the
generators of $G_{n,i}$,
\item the approximate smallest co-volume $V_{n,i}$ of an arithmetic
Kleinian group in which we know $G_{n,i}$ can be embedded.
\end{enumerate}

\bigskip

\begin{center}
\begin{tabular}{|c|c|c|c|c|c|}
\multicolumn{6}{c}{\bf Table 6 - Co-volume of group containing {\boldmath
$G_{3,i}$}} \\ 
\hline
$ i $&$q_{3,i}$&$ d_{kG_{3,i}}$&$ {\rm Ram}_f $&$\delta_{3,i}$&$ V_{3,i} $\\
\hline
$ 1 $&$ z+1 $&$--$&$--$&$ 0 $&$ S_4 $\\
\hline
$ 2 $&$ z^2+3z+1  $&$-- $&$-- $&$ 0 $&$ A_5 $\\
\hline
$ 3 $&$ z^4+6z^3+12z^2+9z+1 $&$ -275 $&$ \emptyset $&$ .1970 $&$ .0390 $\\
\hline
$ 4 $&$ z^4+5z^3+7z^2+3z+1 $&$ -283 $&$ \emptyset $&$ .2108 $&$ .0408 $\\
\hline
$ 5 $&$ z^3+4z^2+4z+2 $&$ -44 $&$ {\cal P}_2 $&$ .2337 $&$ .0661 $\\
\hline
$ 6 $&$ z^3+5z^2+8z+5 $&$ -23 $&$ {\cal P}_5 $&$ .2448 $&$ .0785 $\\
\hline
$ 7 $&$ z^3+3z^2+2z+1 $&$ -23 $&$ {\cal P}_5 $&$ .2480 $&$ .0785 $\\
\hline
$ 8 $&$ z^4+5z^3+8z^2+6z+1 $&$ -563 $&$ \emptyset $&$ .2740 $&$ .1274 $\\
\hline
$ 9 $&$ z^2+3z+3 $&$ -3 $&$ \emptyset $&$.2746 $&$ .0845 $\\
\hline
$ 10 $&$ z^3+4z^2+5z+3 $&$ -31 $&$ {\cal P}_3 $&$ .2770 $&$ .0659 $\\
\hline
$ 11 $&$ z^4+5z^3+6z^2+1 $&$ -491 $&$ \emptyset $&$ .2788 $&$ .1028 $\\
\hline
$ 12 $&$ z^3+4z^2+3z-1 $&$--$&$--$&$ .2831 $&${\rm Fuch.}$\\
\hline
$ 13 $&$ z^6+8z^5+24z^4+35z^3+28z^2+12z+1 $&$ ? $&$ ? $&$ .2944 $&$ ? $\\
\hline
$ 14 $&$ z^3+3z^2+z+1 $&$ -76 $&$ {\cal P}_2 $&$.2970 $&$ .1654 $\\
\hline
\end{tabular}

\bigskip

\begin{tabular}{|c|c|c|c|c|c|}
\multicolumn{6}{c}{\bf Table 7 - Co-volume of group containing {\boldmath 
$G_{4,i}$}} \\
\hline
$ i $&$ q_{4,i}$&$ d_{kG_{4,i}}$&$ {\rm Ram}_f $&$\delta_{4,i}$&$ V_{4,i} $\\
\hline
$ 1 $&$ z+1 $&$-- $&$-- $&$ 0 $&$ S_4 $\\
\hline
$ 2 $&$ z^2+z+1  $&$ -3 $&$ \{{\cal P}_4,{\cal P}_3\} $&$ .4157 $&$ .1268 $\\
\hline
$ 3 $&$ z^3+2z^2+z+1 $&$ -23 $&$ {\cal P}_8 $&$ .4269 $&$ .1374 $\\
\hline
$ 4 $&$ z^2+2z+2 $&$ -4 $&$ \emptyset $&$ .4406 $&$ .2289 $\\
\hline
$ 5 $&$ z^3+3z^2+4z+3 $&$ -31 $&$ {\cal P}_3 $&$ .5049 $&$ .2968 $\\
\hline
$ 6 $&$ z^3+z^2+1 $&$ -31 $&$ {\cal P}_3 $&$ .5225 $&$ .2968 $\\
\hline
$ 7 $&$ z^3+2z^2+2z+2 $&$ -44 $&$ {\cal P}_2 $&$ .5297 $&$ .0661 $\\
\hline
$ 8 $&$  z^3+z^2-z+1 $&$ -44 $&$ {\cal P}_2 $&$ .5297 $&$ .0661 $\\
\hline
$ 9 $&$ z^2+1 $&$ -4 $&$ \emptyset $&$ .5306 $&$ .1526 $\\
\hline
$ 10 $&$ z^2+z-1 $&$--$&$-- $&$ .5306 $&${\rm Fuch.}$\\
\hline
$ 11 $&$ z^4+4z^3+7z^2+6z+1 $&$ -400 $&$ \emptyset $&$ .5306 $&$ .0717 $\\
\hline
$ 12 $&$ z^4+3z^3+4z^2+4z+1 $&$ -331 $&$ \emptyset $&$ .5345 $&$ .4475 $\\
\hline
$ 13 $&$ z^4+z^3-2z^2+1 $&$ -283 $&$ \emptyset $&$ .6130 $&$ .3475 $\\
\hline
\end{tabular}

\bigskip

\begin{tabular}{|c|c|c|c|c|c|}
\multicolumn{6}{c}{\bf Table 8 - Co-volume of group containing {\boldmath 
$G_{5,i}$}} \\
\hline
$ i $&$ q_{5,i}$&$ d_{kG_{5,i}}$&$ {\rm Ram}_f $&$\delta_{5,i}$&$ V_{5,i} $\\
\hline
$ 1 $&$ z^2+3z+1 $&$--$&$-- $&$ 0 $&$A_5$\\
\hline
$ 2 $&$ z^4+5z^3+7z^2+5z+1 $&$ -475 $&$ \emptyset $&$ .4568 $&$ .0933 $\\
\hline
$ 3 $&$ z^4+4z^3+2z^2+z+1 $&$ -275 $&$ \emptyset $&$ .5306 $&$ .0390 $\\
\hline
$ 4 $&$ z^4+3z^3+3z^2+3z+1 $&$ -275 $&$ \emptyset $&$ .6097 $&$ .0390 $\\
\hline
$ 5 $&$ z^2+z-1 $&$-- $&$-- $&$ .6268 $&${\rm Fuch.}$\\
\hline
$ 6 $&$ z^6+6z^5+11z^4+8z^3+7z^2+7z+1 $&$ ? $&$ ? $&$ .6514 $&$ ? $\\
\hline
$ 7 $&$ z^4+5z^3+9z^2+10z+4 $&$ -775 $&$ \{{\cal P}_4,{\cal P}'_4\} $&$
.6717 $&$
.461 $\\ \hline
$ 8 $&$ z^4+6z^3+12z^2+14z+11 $&$ -400 $&$ \emptyset $&$ .6949 $&$ .0717 $\\
\hline
$ 9 $&$ z^4+4z^3+4z^2+5z+5 $&$ -475 $&$ \emptyset $&$ .7195 $&$ .0933 $\\
\hline
$ 10 $&$ z^6+6z^5+12z^4+16z^3+17z^2+8z+1 $&$ ? $&$ ? $&$ .7273 $&$ ? $\\
\hline
$ 11 $&$ z^4+2z^3-2z^2+2z+1 $&$ -400 $&$ \emptyset $&$ .7323 $&$ .0717 $\\
\hline
$ 12 $&$ z^6+5z^5+3z^4-8z^3+z^2+8z+1 $&$ ? $&$ ? $&$ .7725 $&$ ? $\\
\hline
\end{tabular}

\bigskip

\begin{tabular}{|c|c|c|c|c|c|}
\multicolumn{6}{c}{\bf Table 9 - Co-volume of group containing {\boldmath 
$G_{6,i}$}} \\
\hline
$ i $&$ q_{6,i}$&$ d_{kG_{6,i}}$&$ {\rm Ram}_f $&$\delta_{6,i}$&$ V_{6,i} $\\
\hline
$ 1 $&$ z^2+z+1 $&$ -3 $&$ \emptyset $&$ .6584 $&$ .0845 $\\
\hline
$ 2 $&$ z^2+1 $&$ -4 $&$ \{{\cal P}_9,{\cal P}_2\} $&$ .7642 $&$ .3053 $\\
\hline
$ 3 $&$ z^2-z+1 $&$ -3 $&$ \emptyset $&$ .8314 $&$ .1691 $\\
\hline
$ 4 $&$ z^2+z+2 $&$ -7 $&$ \{{\cal P}_2,{\cal P}'_2\} $&$ .8500 $&$ .5555 $\\
\hline
$ 5 $&$ z^3+z^2+2z+1 $&$ -23 $&$ {\cal P}_{27} $&$ .8539 $&$ .5106 $\\
\hline
$ 6 $&$ z^3+z+1 $&$ -31 $&$ {\cal P}_3 $&$ .8786 $&$ .3298 $\\
\hline
$ 7 $&$ z-1  $&$--$&$-- $&$ .8813 $&$\rm Fuch.$\\
\hline
$ 8 $&$ z^3-z^2+1 $&$ -23 $&$ {\cal P}_{27} $&$ .9106 $&$ .5106 $\\
\hline
\end{tabular}

\bigskip

\begin{tabular}{|c|c|c|c|c|c|}
\multicolumn{6}{c}{\bf Table 10 - Co-volume of group containing {\boldmath 
$G_{7,i}$}} \\
\hline
$ i $&$ q_{7,i}$&$ d_{kG_{7,i}}$&${\rm Ram}_f $&$\delta_{7,i}$&$ V_{7,i} $\\
\hline
$ 1 $&$ z^3+4z^2+3z-1 $&$-- $&$-- $&$ .5452 $&${\rm Fuch.}$\\
\hline
$ 2 $&$ z^6+7z^5+17z^4+21z^3+17z^2+7z+1  $&$ ?  $&$ ?  $&$ .8162 $&$ ?  $\\
\hline
$ 3 $&$ z^3+z^2-2z-1 $&$--$&$--$&$ 1.0704 $&${\rm Fuch.}$\\
\hline
\end{tabular}

\end{center}

\bigskip

\centerline{\bf Notes}
\begin{enumerate}
\item In these tables, we have adopted the notation of labeling a prime ideal of norm
$p^n$ by ${\cal P}_{p^n}$. In the cases considered, this uniquely determines the ideal
except for the indicated cases where ${\cal P}_2,{\cal P}_2'$ or ${\cal P}_4,{\cal
P}_4'$ both appear. \item We have not dealt with the cases where 
$[{\bf Q}(\gamma) : {\bf Q}] > 4$.
\end{enumerate}
\bigskip

\centerline{\bf Remarks}
\begin{enumerate}
\item As remarked in \S 7, any subgroup of an arithmetic Kleinian group is contained in a
group $\Gamma_{S,S'}$ \cite{Bo}. In all cases when $n = 3$, applying this to $G_{3,i}$
we have $S = \phi$. That is not the situation in many cases when $n = 4,5,6$, but
in most cases $S$ consists of a single prime ideal. The volume of $\Gamma_{S,S'}$ is
readily calculated from the formulas given in \cite{Bo}.

\item A detailed knowledge of subgroups commensurable with tetrahedral groups enable
some groups to be identified immediately by their commutator parameter  and so the
minimum volume of a group in which they lie is quickly determined. (See comments after
the remarks on the tables). The arithmetic nature of tetrahedral groups is discussed in
\cite{MR2} and \cite{R1} from which much of the algebraic data of Table 8 is taken.

\item {$\bf \hbox {Table 7, n=4.}$} 
The order we have considered is (suppressing subscripts) 
\[ R_{kG}[1,f^2,(gf)^2,f^2(gf)^2] \]
which has discriminant $\langle 2\gamma(\gamma+2) \rangle$.  
Only in the cases $i = 2,3$ is this order maximal. In these cases, the group
$G_{4,i}$ is contained in a $\Gamma_{\phi,\phi}$ as is also the case for 
$i = 7,8$ following an argument similar to that used in the case of $G_{3,10}$ in \S 7.
In the cases $i = 4,5,6,12,13$, $S$ consists of a single prime, while in the 
case $i = 9$, $S$ may consist of two primes. The groups $G_{4,1}$ and $G_{4,11}$ are
subgroups of groups commensurable with tetrahedral groups.

\item {$\bf \hbox {Table 8, n = 5.}$}  Referring to Table 3,  note that the
polynomials $p_{5,i}$ admit $z + 1$ as a factor in the cases i = 3,8,9,10,11,12.
Using this, the field discriminants of Table 8 are readily determined.

The order here is chosen to be 
$$ R_{kG}[1,gfg^{-1},f^{-1},[g,f]]$$
which has discriminant $\langle \gamma(\gamma - \beta) \rangle$ where
$\beta = (-5 + \surd 5)/2$. From this, as before, we determine the finite ramification
of the algebra, which here, is only non-empty in the case $i = 7$. All other groups
arise in considering cocompact tetrahedral groups.

\item {$\bf \hbox {Table 9, n = 6.}$} Now we consider the order 
$$R_{kG}[1,f^2,gf^{-2}g^{-1},f^2gf^{-2}g^{-1}]$$
which has discriminant $\langle 9\gamma(\gamma + 1) \rangle$. For $i = 2,5,8$, the order
is maximal so that $S = \phi$, while when $i = 3,4,6$, $S$ consists  of at most one
prime ideal. The group $G_{6,1}$ is commensurable with a tetrahedral group.  See also
\cite{GM2}.  
\end{enumerate}

\bigskip

We now comment on the relationship between some of these groups and tetrahedral groups.
Recall by a tetrahedral group we mean the orientation-preserving subgroup of index $2$
in the group generated by reflections in the faces of tetrahedron in ${\bf H}^3$, where
some vertices may be ideal. Justifications for the comments below can be deduced for
example from \cite{R1} where, as here, the notation for the tetrahedra is that of
\cite{Be}.

The groups $\PGL(2,{\bf Z}[i])$ and $\PGL(2,O_3)$ are tetrahedral, the tetrahedrons
being $T[3,2,2;4,2,4]$ and $T[3,2,2;6,2,3]$. In the tables these arise in the cases of
$G_{4,9}$ and $\PGL(2,O_3)$
occurs in the cases $G_{3,9}$ and $G_{6,1}$.
Other non-cocompact groups which arise although not tetrahedral, contain a
tetrahedral group of index 2. The case $G_{4,4}$ yields a group containing 
the tetrahedral group $T[4,2,2;4,2,4]$ of index $2$, and the case $G_{6,3}$ yields
a group containing  the tetrahedral group $T[2,2,3;2,6,3]$ of index $2$.

In the case of order $5$ the algebras defined over the quartic fields listed 
in Table 8 except $i=7$ yield groups 
commensurable with cocompact tetrahedral groups. In the 
case where the discriminant is $-400$ the tetrahedral group is $T[2,2,4;2,3,5]$
(which  also arises in the case $G_{4,11}$), the discriminants $-275$ and $-475$
yield groups commensurable with the tetrahedral groups $T[2,2,3;2,5,3]$ and
$T[2,2,5;2,3,5]$ respectively. The case $T[2,2,3;2,5,3]$ also arises in the case
$G_{3,3}$.

Finally the order $6$ case $G_{6,4}$
yields a group
commensurable with the tetrahedral group $T[2,3,4;2,3,4]$.

\bigskip
Recently K. N. Jones and the fourth author have developed a computer program to
study explicitly how the geometry and topology of certain arithmetic Kleinian
groups varies as the number theoretic data is varied. Among other things it constructs a
fundamental polyhedron for the action of certain unit groups of orders and their
normalizers in ${\bf H}^3$, it determines the co-volume and it computes a presentation.
We present in Table 11 the results produced by that program for the two-generator groups
discussed  in this paper.  Since all but $G_{4,11}$, $G_{5,7}$, $G_{6,4}$ and the
Fuchsian examples turn out to have finite co-volume, we conclude that, with these
exceptions, all these  groups are themselves arithmetic.

\begin{theorem}
The following two-generator Kleinian groups $G_{n,i}$ generated by elliptics of
order $n$ and $2$ are arithmetic.
\begin{enumerate}
\item $G_{3,i}$,  $i=3,\ldots, 14$, $i\neq 12, 13$,
\item $G_{4,i}$,	 $i=2,\ldots, 13$, $i\neq 10, 11$,
\item $G_{5,i}$,	 $i=2,\ldots, 11$, $i\neq 5,6,7,10$,
\item $G_{6,i}$,	 $i=1,\ldots, 8$, $i\neq 4$.
\end{enumerate}
\end{theorem}

\bigskip

\begin{center}
\begin{tabular}{|c|c|c|c|c|}
\multicolumn{5}{c}{\bf Table 11 - Co-volumes of {\boldmath $G_{n,i}$}} \\
\hline
$ i $&$ G_{3,i}$&$ G_{4,i}$&$ G_{5,i} $&$G_{6,i}$\\
\hline
$ 1 $&$  S_4 $&$ S_4 $&$ A_5 $&$ .2537 $\\
\hline
$ 2 $&$ A_4   $&$ .2537 $&$ .0933 $&$ .6106 $\\
\hline
$ 3 $&$ .0390 $&$ .1374 $&$ .0390 $&$ .5074 $\\
\hline
$ 4 $&$ .0408 $&$ .4579 $&$ .4686 $&$ \infty $\\
\hline
$ 5 $&$ .1323 $&$ .5936 $&$ {\rm Fuch.} $&$ 1.0212 $\\
\hline
$ 6 $&$ .1571 $&$ .5936 $&$ ?     $&$ 1.3193 $\\
\hline
$ 7 $&$ .1571 $&$ .7943 $&$ \infty $&${\rm Fuch.} $\\
\hline
$ 8 $&$ .1274 $&$ .2647 $&$ .8612 $&$ 1.0212 $\\
\hline
$ 9 $&$ .3383 $&$ .9159 $&$ 1.1199 $&$ -- $\\
\hline
$ 10 $&$ .2638 $&$ {\rm Fuch.} $&$ ? $&$ -- $\\
\hline
$ 11 $&$ .2056 $&$ \infty $&$ .8612 $&$ -- $\\
\hline
$ 12 $&${\rm Fuch.} $&$ .8951 $&$ ? $&$ -- $\\
\hline
$ 13 $&$ ?  $&$ .3475 $&$ -- $&$ -- $\\
\hline
$ 14 $&$ .3308 $&$ -- $&$ -- $&$ -- $\\
\hline
\end{tabular}
\end{center}

\bigskip
We conclude this section with some philosophical remarks.  Several of the algebras and
small volume groups arise in the considerations of \cite{CFJR}. This work identifies
the smallest volume of an orientable arithmetic hyperbolic 3-manifold, namely that of
the Weeks manifold. The content of this article together with \cite{CF} and the works
\cite{GM1}--\cite{GM8} seem to indicate a convergence of arithmetic and geometric
ideas. In particular the smallest volume orientable arithmetic orbifold arises as the
quotient of hyperbolic 3-space  ${\bf H}^3$ by a ${\bf Z}_2$-extension of the
tetrahedral group $T[2,3,3;2,5,3]$ (cf. \cite{CF}) and it would seem that this
article and \cite{GM1}---\cite{GM6} are ``converging on'' the arithmetic orbifold
just described as the overall smallest one, as is conjectured.

\section{Criteria for simple axes}

Suppose that $f$ is an elliptic element of order $n$ in a Kleinian group $G$.  
Then $f$ is {\it simple} if axis($f$) is precisely invariant, that is, if for
each $h \in G$  
\[ h({\rm axis}(f)) = {\rm axis}(hfh^{-1})\]
either coincides with or is disjoint from axis($f$). If $f$ is not simple, then
one of the following is true \cite{Bea} and \cite{Mas}.  
\begin{enumerate}
\item The axes of $f$ and $hfh^{-1}$ intersect in ${\bf H}^3$ and $G$ contains a
subgroup isomorphic to the regular solid group $A_4$, $S_4$ or $A_5$.  
\item The axes of $f$ and $hfh^{-1}$ have one endpoint in common and $G$ contains a
subgroup isomorphic to the (3,3,3), (2,4,4) or (2,3,6) euclidean triangle group.  
\end{enumerate}

We show here that if $G$ is a Kleinian group of finite co-volume, then each of
these conditions forces the invariant quaternion algebra of $G$ to satisfy certain
criteria. Hence if the algebra fails to satisfy these criteria, then every 
elliptic in $G$ of order $n \geq 3$ must be simple.  

\bigskip

\centerline{\bf Case where axes of {\boldmath $f$} and {\boldmath $hfh^{-1}$}
intersect in ${\bf H}^3$}

\medskip
Let ${\cal H}$ denote Hamilton's quaternions so that
${\cal H} =  \displaystyle{\biggl ({-1,-1 \over {\bf R}}\biggr)}$ and let $\sigma$ 
denote the embedding $\sigma : {\cal H}^1 \rightarrow \SL(2,{\bf C})$  given by 
\[\sigma(a_0 + a_1i + a_2j + a_3ij) = \pmatrix{a_0 + a_1i & a_2 + a_3i \cr 
-a_2 + a_3i & a_0 - a_1i},\]
where ${\cal H}^1$ is the group of elements of norm 1. If $n$ denotes the norm on
${\cal H}$, then there is an epimorphism 
\[\Phi : {\cal H}^1 \rightarrow \SO(3, {\bf R})\] 
where $\SO(3,{\bf R})$ is represented as the orthogonal group of the quadratic
subspace $V$ of ${\cal H}$ spanned by $\{i,j,ij\}$, that is, the  pure quaternions, 
equipped with the restriction of the norm form, so that 
$n(x_1i + x_2j + x_3ij) = x_1^2 + x_2^2 + x_3^2.$ The mapping $\Phi$ is defined by
$\Phi(\alpha) = \phi_{\alpha}$ where  
\[\phi_{\alpha}(\beta) = \alpha \beta
\alpha^{-1} \quad \alpha \in {\cal H}^1  \quad \beta \in V.\]
The kernel of $\Phi$ is $\{ \pm 1\}$. The binary tetrahedral group is a 
central extension of an element of order 2 by the tetrahedral group
and can be faithfully represented in ${\cal H}$. This is also true 
for the binary octahedral group and the binary icosahedral group as will 
now be shown.

If the tetrahedron has its vertices at 
\[i + j + ij, \;\;\; i - j - ij, \;\;\; -i + j - ij, \;\;\; -i - j + ij,\] 
then $\phi_{\alpha_1}$ is a rotation of order 2 about the axis  through the edge
mid-point $i$ if $\alpha_1 = i$, while $\phi_{\alpha_2}$  is a rotation of order 3
about the axis through the vertex $i + j + ij$  when $\alpha_2 = (1 + i + j +
ij)/2$. The binary tetrahedral group $\Gamma_1$  is thus generated by $\alpha_0 =
-1,\alpha_1,\alpha_2$ in ${\cal H}^1$.

If the cube has its vertices at
\[\pm i \pm j \pm ij,\] 
then  $\phi_{\alpha_3}$ is a rotation of order 4 about the axis through the
mid-point  $i$ of a face when $\alpha_3 = (1 + i)/\sqrt{2}$. Thus the binary
octahedral  group $\Gamma_2$ is generated  by $\alpha_0,\alpha_2,\alpha_3$.

If the regular icosahedron has its vertices at 
\[\pm i \pm j \pm ij, \;\;\; \pm \tau i \pm \tau^{-1} j, \;\;\; \pm \tau j \pm \tau^{-1}
ij, \;\;\;  \pm \tau ij \pm \tau^{-1}i,\] 
then $\phi_{\alpha_4}$ is a rotation of order 5 about the axis through the
mid-point of the face with  vertices 
\[i + i + ij, \;\;\; i + j - ij, \;\;\; \tau i + \tau^{-1}j, \;\;\; \tau j +
\tau^{-1}ij, \;\;\; \tau j - \tau^{-1} ij,\]  
where $\alpha_4 = (\tau + \tau^{-1}i + j)/2$ and $\tau =
(\sqrt{5} + 1)/2$. The binary icosahedral group $\Gamma_3$ is then  generated by
$\alpha_0,\alpha_2,\alpha_4$.

The groups $P\sigma(\Gamma_1) \cong A_4, P\sigma(\Gamma_2) \cong S_4, 
P\sigma(\Gamma_3) \cong A_5$, where $P$ is the projection 
$P : SL(2,{\bf C}) \rightarrow PSL(2,{\bf C})$,  are said to be in standard
form in  $PSL(2,{\bf C})$. If $G$ is a discrete group 
which contains a finite subgroup $F$ isomorphic to one of these regular solid 
groups, then $G$ can be conjugated so that $F$ is in {\it standard form}.
Note that if $G$ contains any of these finite groups, it will contain 
a subgroup isomorphic to $A_4$.

\begin{lemma} 
Let $G$ be a Kleinian group of finite co-volume with invariant quaternion algebra
$A$ and number field  $k$. If $G$ contains a subgroup isomorphic to $A_4$, 
 then
\begin{equation}   
A \cong \displaystyle{\biggl({-1,-1 \over k}\biggr)}.
\end{equation}
In particular, the only finite primes at which $A$ can be ramified are the dyadic
primes. 
\end{lemma}
{\bf Proof.} 
Suppose that $G$ contains a subgroup isomorphic to $A_4$. Since $A_4$ has no
subgroup of index $2$ and $A_5$ is simple, we have $A_n = A_n^{(2)}$ for $n=4$ or
$5$. So if $A_n \subset G$, then  $A_n \subset G^{(2)}$.  Thus by conjugation, we
can assume that $\sigma(\Gamma_1) \subset {\cal G}$ where  $P{\cal G}= G^{(2)}$.
Now   \[A = \{ \sum a_i g_i: a_i \in k, \quad g_i \in  {\cal G} \}.\] 
Let $\displaystyle{A_0 = \biggl({-1,-1 \over {\bf Q}}\biggr)}$.  Then 
\[A_0 \cong \{ \sum a_i g_i : a_i \in {\bf Q} \quad,  g_i \in
\sigma(\Gamma_1) \}\] since $1,i,j,ij \in \Gamma_1$. Now the quaternion algebra 
\[\{ \sum a_i g_i : a_i \in k \quad g_i \in  \sigma(\Gamma_1)\}\] 
lies in $A$, is isomorphic to  $\displaystyle{A_0 \otimes_{\bf Q} k}$ and is
4-dimensional. Thus   
\[A \cong \biggl({-1,-1 \over k}\biggr).\]

Finally it follows from this form of the Hilbert symbol, that $A$ splits over all ${\cal P}$-adic fields $k_{\cal P}$, 
with ${\cal P}$ non-dyadic, that is, primes not dividing $2$. $\Box$

\begin{lemma}
Let $G$ be as in {\rm Lemma 9.1}. If $G$ contains a subgroup isomorphic to
$A_5$ and if $[k : {\bf Q}] = 4$, then $A$ has no finite ramification.
\end{lemma}

{\bf Proof.} $A$ can at worst have dyadic finite ramification. 
Note that, since $G$ must contain an element of order 5, ${\bf Q}(\surd 5) 
\subset k$. There is a unique prime ${\cal P}$ in ${\bf Q} (\sqrt{5})$ such that ${\cal
P} \mid 2$. So if ${\cal P}$  ramifies or is inert in $k \mid {\bf Q}(\sqrt{5})$, then
there will only be one dyadic prime in $k$ at which $A$ cannot be  ramified for parity
reasons. Suppose then that ${\cal P}$  splits as ${\cal P}_1{\cal P}_2$ so that
$k_{{\cal P}_1} \cong  k_{{\cal P}_2} \cong {\bf Q}(\sqrt{5})_{{\cal P}}$. But,
again for parity reasons, the quaternion algebra  
\[\biggl(\frac{-1,-1}{{\bf Q}(\sqrt{5})}\biggr)\]  
splits in  the field ${\bf Q}(\sqrt{5})_{\cal P}$. Hence  
\[\biggl(\frac{-1,-1}{k}\biggr) \]
splits in $k_{{\cal P}_1}$ and $k_{{\cal P}_2}$ and $A$ has no finite
ramification. $\Box$ 

\bigskip

Lemma 9.1 has the following partial converse.
\begin{lemma} 
Suppose that $k$ has exactly one complex place. If  
\[{A \cong \biggl(\frac{-1,-1}{k}\biggr)},\] 
then there is an arithmetic Kleinian group in the commensurability class defined by
$A$  which contains $S_4$. Furthermore if ${\bf Q}(\sqrt{5})
\subset k$, then there is a group in the commensurability class containing $A_5$.  
\end{lemma}
{\bf Proof.} To exhibit a group containing $S_4$ we proceed as follows. Let 
$R_k$ denote the ring of integers of $k$ and let 
\[{\cal O} = R_k[1,i,j,(1 + i + j + ij)/2].\]  
It is easily checked that ${\cal O}$ is an order in $A$. Note that $\Gamma_1
\subset  {\cal O}^1$. Furthermore the element  $1 + i \in {\cal O}$  normalizes
${\cal O}^1$. Thus if $\rho$ is a representation of  $A$ into $M(2,{\bf C})$, then  
\[P\rho({\cal O}^1) \subset \langle P\rho ({\cal O}^1), \;\;\; 
P\rho (1 + i)\rangle = G\]  
as a subgroup of index 2. But $P\rho(1 + i)$ also normalizes
$P\rho(\Gamma_1)$ and  \[\langle P\rho(\Gamma_1) , P\rho(1 + i)\rangle \cong S_4\]
from the description above. Thus $S_4 \subset G$.

Next, in the case where ${\bf Q}(\surd 5) \subset k$, if we let 
\[{\cal O} = R_k[1,i,1/2(\tau + \tau^{-1}i + j),1/2(-\tau^{-1} 
+ \tau i + ij)],\] 
then ${\cal O}$ is an order in $A$. Note that 
the sum of the last two elements above is $\tau^{-1} + 
1/2(1 + i + j + ij)$ so that $\Gamma_3 \subset {\cal O}^1$. Hence 
$A_5 \subset P\rho({\cal O}^1)$ as required. $\Box$

\bigskip

\centerline{\bf Case where axes of {\boldmath $f$} and {\boldmath $hfh^{-1}$} have
one common endpoint}

\medskip

This case is handled by the following result.
\begin{lemma} 
Let $G$ be a subgroup of an arithmetic Kleinian group with invariant quaternion
algebra $A$ over $k$. If $G$ contains a subgroup isomorphic to the {\rm(3,3,3)}, 
{\rm (2,4,4)} or {\rm (2,3,6)} euclidean triangle group, then
\begin{equation} 
A \cong M(2,{\bf Q}(\sqrt{-1})) \; or \; A \cong M(2,{\bf Q}(\sqrt{-3})).
\end{equation}
\end{lemma}
{\bf Proof.}
Since $G$ contains a parabolic, a well known criterion for non-cocompact 
arithmetic groups implies that $A \cong M(2,{\bf}(\sqrt{-d}))$
for some square-free integer $d$.  Then it follows easily from \cite{NR} or 
\cite{R3} that ${\bf Q}(\sqrt{-1})$ or ${\bf Q}(\sqrt{-3})$ is contained in $kG$
and we obtain (9.6). $\Box$ 

\bigskip
Lemmas 9.1, 9.3 and 9.5 yield necessary conditions for the existence of a
non-simple elliptic.  We conclude this section with a sufficient condition
for the existence of a non-simple elliptic. 

\begin{lemma}
Suppose that $f$ and $h$ are M\"{o}bius transformations and that $f$ is elliptic of
order $n \geq 3$. If 
\[\gamma(f,h) \in (\beta(f),0)\]
or if
\[\gamma(f,h) = \beta(f) \,\,\, and \,\,\, \beta(h) \neq -4,\]
then the axes of $f$ and $hfh^{-1}$ intersect in a single point.
\end{lemma} 
{\bf Proof.}  Let $\beta=\beta(f)$ and $\gamma=\gamma(f,h)$.
Then 
\begin{equation}
\gamma(f,hfh^{-1}) = \gamma(\gamma -\beta) 
\end{equation} 
and by Lemma 2.1,
\begin{eqnarray*}
2 \cosh(\delta(f,hfh^{-1}))^2 & = & \cosh(2\delta(f,hfh^{-1})) + 1 \\
       & = &\left|\frac{4\gamma(f,hfh^{-1})}{\beta(f) \beta(hfh^{-1})}+1\right|+
            \left|\frac{4\gamma(f,hfh^{-1})}{\beta(f)\beta(hfh^{-1})}\right| +1 \\
       & = &\left|\frac{4\gamma(\gamma-\beta)+\beta^2}{\beta^2}\right|+
               \left|\frac{4\gamma(\gamma-\beta)}{\beta^2}\right| \\
       & = & 2\left(\frac{|\gamma-\beta|+|\gamma|}{|\beta|}\right)^2
\end{eqnarray*}
from which we obtain
\begin{equation}
\cosh(\delta(f,hfh^{-1})) = \frac{|\gamma-\beta|+|\gamma|}{|\beta|}.
\end{equation}

If $\gamma \in (\beta,0)$, then $f$ and $hfh^{-1}$ have disjoint fixed
points by (9.8), $\delta(f,hfh^{-1})= 0$ by (9.9) and hence {\rm axis($f$) and 
{\rm axis($hfh^{-1}$) intersect in a single point in ${\bf H}^3$.  If $\gamma =
\beta$ and if $\beta(h) \neq -4$, then (9.8) and the fact that $h$ is not of order 2
imply that $f$ and $hfh^{-1}$ have a single fixed point in common. $\Box$

\section{Simple axes in extremal groups}

Finally, motivated by the search for lower bounds for the volume of hyperbolic
3-orbifolds, we apply the criteria of \S 9 to determine whether or not the
generator $f$ of order $n$ in the groups $G_{n,i}$ is simple.  The results
are summarized in Table 12 at the end of this section.

\bigskip
\centerline{\bf Case n=3}
\medskip
From Table 6 in \S 8 and Lemma 9.1 we can deduce that for $i = 6,7,10$ the groups
$G_{3,i}$ do not have subgroups isomorphic to $A_4$, $S_4$ or $A_5$. In addition,
$G_{3,12}$ is Fuchsian. Hence $f$ is simple in these four cases.  For the other cases
we apply Lemma 9.7 with $\beta(f)=-3$ to show that $f$ is not  simple by exhibiting
a suitable element $h$ in $G_{n,i}$.
    
\begin{enumerate}
\item $G_{3,1}$:  $\gamma(f,g)=-1$ and $\langle f,g \rangle$ is isomorphic to
$S_4$.  
\item $G_{3,2}$:  $\gamma(f,g)=-.3819\ldots$ and $\langle f,g \rangle$ is
isomorphic to $A_5$.
\item $G_{3,3}$:  If $h=gfg$, then $\gamma(f,h)=-2.618\ldots$. 
\item $G_{3,4}$:  If $h=gfgfgfg$, then $\gamma(f,h)=-1$. 
\item $G_{3,5}$:  If $h=gfgfg$, then $\gamma(f,h)=-2$.
\item $G_{3,8}$:  If $h=gfgfgf^{-1}gf^{-1}g$, then $\gamma(f,h)=-2$.
\item $G_{3,9}$:  If $h=gfg$, then $\gamma(f,h)=-3$ and $\beta(h)=-3$.
\item $G_{3,11}$:  If $h=gfgfgfgfg$, then $\gamma(f,h)=-1$.
\item $G_{3,13}$:  If $h=kfk^{-1}fk$ where $k=gfgfg$. then $\gamma(f,h)=-1$.
\item $G_{3,14}$:  If $h=gfgfgfg$, then $\gamma(f,h)=-2$.
\end{enumerate}

\bigskip
\centerline{\bf Case n=4}
\medskip

The methods here are the same as for $n = 3$. In particular, the cases where
$i = 2,5,6$ are dealt with by Lemma 9.1 as above while $G_{4,10}$ is Fuchsian. 
Hence $f$ is simple in these four cases.  Otherwise Lemma 9.7 with $\beta(f)=-2$
shows that $f$ is not simple when $i \neq 7,11$.  

\begin{enumerate}
\item $G_{4,1}$: $\gamma(f,g)=-1$ and $\langle f,g \rangle$ is isomorphic to $S_4$.
\item $G_{4,3}$: If $h=gfgfg$, then $\gamma(f,h) = -1$.
\item $G_{4,4}$: If $h=gfg$, then $\gamma(f,h) = -2$ and $\beta(h)=-2$.
\item $G_{4,8}$: If $h=gfgfgfg$, then $\gamma(f,h) = -1$.
\item $G_{4,9}$: If $h = gfgfg$, then $\gamma(f,h) = -2$ and $\beta(h)=2i$.
\item $G_{4,12}$: If $h=gfgfgf^{-1}gf^{-1}g$, then $\gamma(f,h) = -1$.
\item $G_{4,13}$: If $h=gfgfgfgfg$, then $\gamma(f,h) = -1$.
\end{enumerate}

\bigskip 
\centerline{\bf Case n=5}
\medskip

In this case, the only group which can be dealt with directly by the methods of 
\S 9 is $G_{5,7}$. In this example, although there is only dyadic ramification, the
algebra is not isomorphic to  
\[\displaystyle \biggl({-1,-1 \over k_{5,7} } \biggr)\] 
as is shown by  Lemma 9.3.

Finally since $G_{5,5}$ is Fuchsian we conclude that $f$ is simple when $i=4,5,7$. 
We can then apply Lemma 9.7 with  $\beta(f)=-1.3819\ldots$ to show that $f$
is not simple when $i \neq 8,9,11$.

\begin{enumerate}
\item $G_{5,1}$: $\gamma(f,g)=-.38196\ldots$ and $\langle f,g \rangle$ is isomorphic
to $A_5$. 
\item $G_{5,2}$: If $h=gfg$, then $\gamma(f,h) = -1$.
\item $G_{5,3}$: If $h=gfgfg$, then $\gamma(f,h) = -.38196\ldots$.
\item $G_{5,6}$: If $h=gfgfg$, then $\gamma(f,h) = -1$.
\item $G_{5,10}$: If $h=gfgfgf^{-1}gf^{-1}g$, then $\gamma(f,h) = -.38196\ldots$.
\item $G_{5,12}$: If $h=gfgfgfg$, then $\gamma(f,h) = -1$.
\end{enumerate}

\bigskip
\centerline{\bf Case n=6}
\medskip

When $n = 6$, Lemma 9.5 implies that the only cases for which $f$ cannot be
simple are $G_{6,1}$  and $G_{6,3}$.  The calculations for these using Lemma 9.7
with $\beta(f)=-1$ are given below.   
\begin{enumerate}
\item $G_{6,1}$: If $h=gfg$, then $\gamma(f,h)=-1$ and $\beta(h)=-1$.
\item $G_{6,3}$: If $h=gfgfg$, then $\gamma(f,h) = -1$ and $\beta(h)=.5+i2.598\ldots$.
\end{enumerate}

\bigskip

The six groups $G_{4,7}$, $G_{4,11}$, $G_{5,4}$, $G_{5,8}$, $G_{5,9}$ and $G_{5,11}$ 
have not been treated in the above calculations.  The program of Jones and Reid
mentioned in \S 8 shows that $f$ is simple in all of these groups.

\bigskip

We summarize these results in Table 12 below.
\begin{center}
\begin{tabular}{|c|c|c|c|c|}
\multicolumn{5}{c}{\bf Table 12 - {\boldmath $f$} simple elliptic} \\
\hline
$ i $&$ G_{3,i}     $&$ G_{4,i}     $&$ G_{5,i}     $&$ G_{6,i}    $\\
\hline
$ 1 $&$ S_4         $&$ S_4         $&$ A_5         $&$ {\rm No}   $\\
\hline
$ 2 $&$ A_4         $&$ {\rm Yes}   $&$ {\rm No}    $&$ {\rm Yes}  $\\
\hline 
$ 3 $&$ {\rm No}    $&$ {\rm No}    $&$ {\rm No}    $&$ {\rm No}   $\\
\hline
$ 4 $&$ {\rm No}    $&$ {\rm No}    $&$ {\rm Yes}   $&$ {\rm Yes}  $\\
\hline
$ 5 $&$ {\rm No}    $&$ {\rm Yes}   $&$ {\rm Fuch.} $&$ {\rm Yes}  $\\
\hline
$ 6 $&$ {\rm Yes}   $&$ {\rm Yes}   $&$ {\rm No}    $&$ {\rm Yes}  $\\
\hline
$ 7 $&$ {\rm Yes}   $&$ {\rm Yes}   $&$ {\rm Yes}   $&$ {\rm Fuch.}$\\
\hline
$ 8 $&$ {\rm No}    $&$ {\rm No}    $&$ {\rm Yes}   $&$ {\rm Yes}  $\\
\hline
$ 9 $&$ {\rm No}    $&$ {\rm No}    $&$ {\rm Yes}   $&$ --         $\\
\hline
$ 10 $&$ {\rm Yes}  $&$ {\rm Fuch.} $&$ {\rm No}    $&$ --         $\\
\hline
$ 11 $&$ {\rm No}   $&$ {\rm Yes}   $&$ {\rm Yes}   $&$ --         $\\
\hline
$ 12 $&${\rm Fuch.} $&$ {\rm No}    $&$ {\rm No}    $&$ --         $\\
\hline
$ 13 $&$ {\rm No}   $&$ {\rm No}    $&$ --          $&$ --         $\\
\hline
$ 14 $&$ {\rm No}   $&$ --          $&$ --          $&$ --         $\\
\hline
\end{tabular}
\end{center}

\section{\bf Commutator parameter diagrams}

We conclude this paper with four diagrams which result from the disk covering
argument described in \S 2 for the cases where $n=3,4,5,6$.  They represent the
{\em only} possible values for the commutator parameter $\gamma=\gamma(f,g)$ of a
discrete group $\langle f,g \rangle$ for which

1. $f$ is an elliptic of order $n$,

2. $\gamma$ lies in the union of the indicated disks. 

\noindent The values in these diagrams which correspond to  subgroups of arithmetic
groups are listed in Tables 1, 2, 3 and 4 in \S 6.  As noted earlier in \S 2, {\em all}
of the values in the $n=3$ diagram have this property.

\centerline{\bf Only possible values for commutator parameter when n=3}
\nopagebreak
\psfig{file=Figure3.eps,width=5.42in}

\bigskip

\centerline{\bf Only possible values for commutator parameter when n=4}
\nopagebreak
\psfig{file=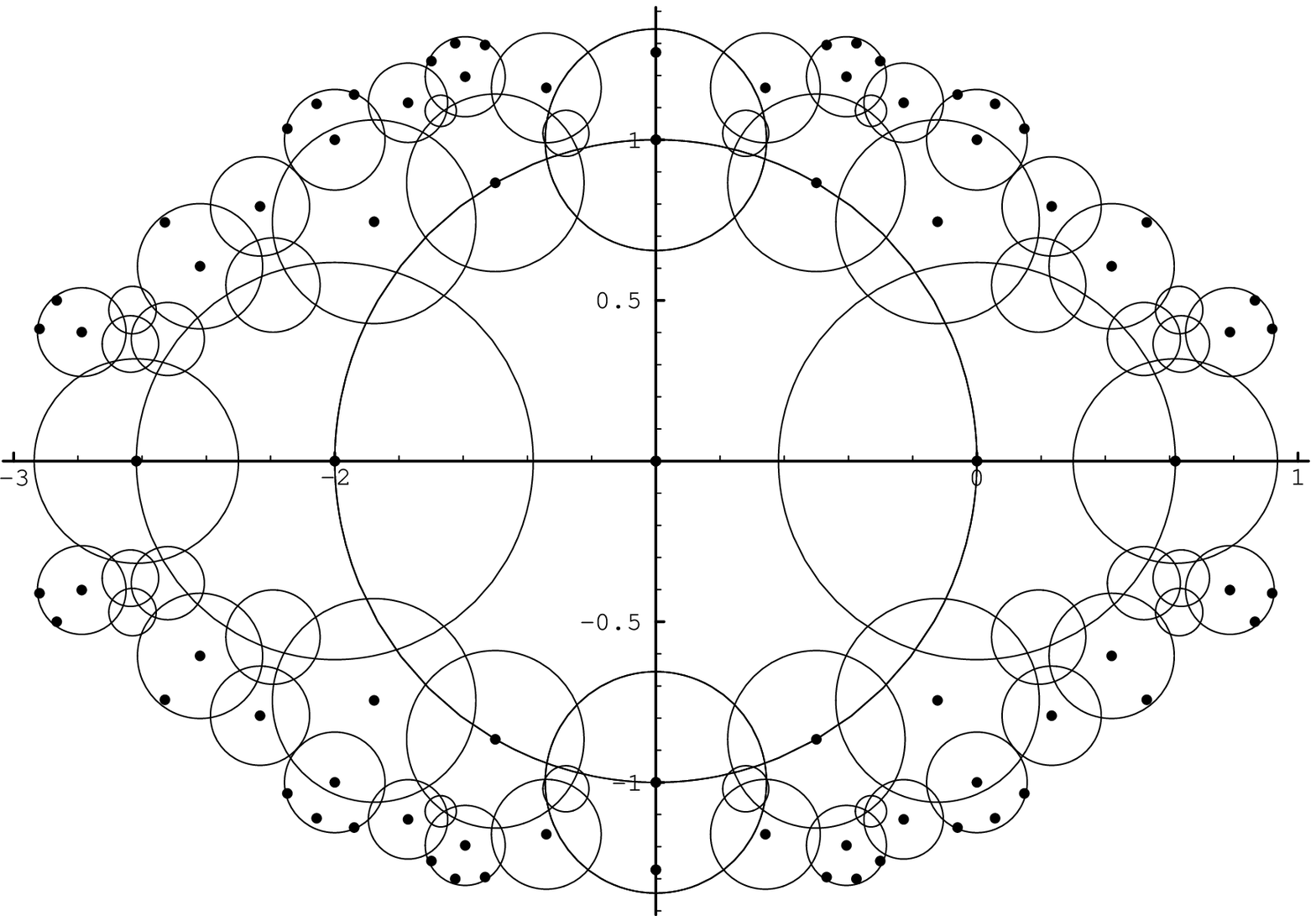,width=5.42in}

\bigskip

\centerline{\bf Only possible values for commutator parameter when n=5}
\nopagebreak
\psfig{file=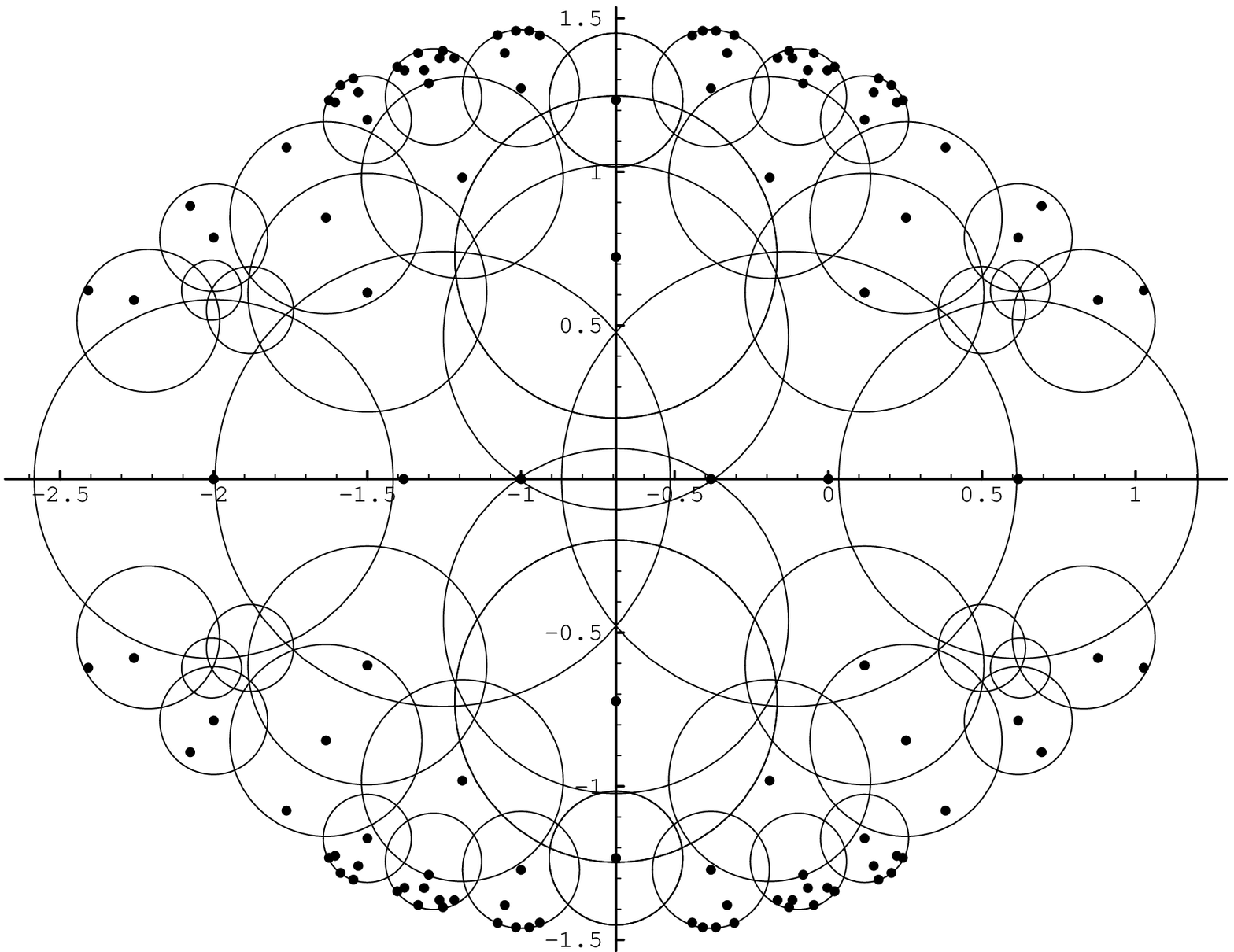,width=5.42in}

\bigskip

\centerline{\bf Only possible values for commutator parameter when n=6}
\nopagebreak
\psfig{file=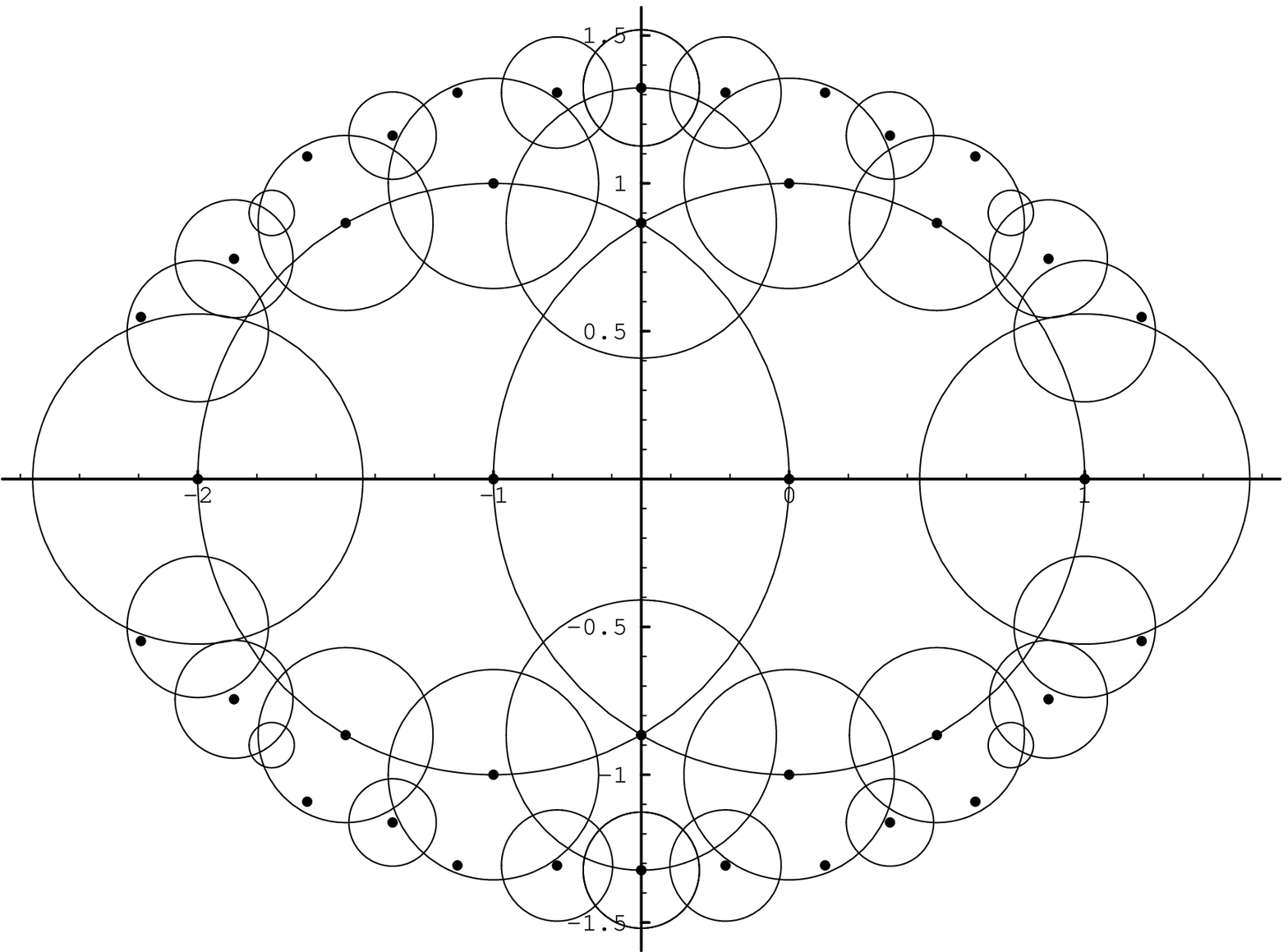,width=5.42in}

\centerline{\bf Authors' addresses}

\noindent F. W. Gehring 

University of Michigan, Ann Arbor, Michigan, USA

\medskip
\noindent C. Maclachlan

University of Aberdeen, Aberdeen, Scotland

\medskip
\noindent G. J. Martin

University of Auckland, Auckland, New Zealand

Australian National University, Canberra, Australia

\medskip
\noindent A. W. Reid

University of Cambridge, Cambridge, England

\end{document}